\documentclass{amsproc}
\usepackage{bbm}

\usepackage{amsfonts,amsmath,amsthm,amssymb,latexsym,mathrsfs,stmaryrd}
\usepackage{hyperref}
\usepackage{mathtools}
\usepackage{graphicx}
\usepackage{subcaption}
\usepackage[bb=libus]{mathalpha}

\makeatletter
\newcommand{\xdashrightarrow}[2][]{\ext@arrow 0359\rightarrowfill@@{#1}{#2}}
\usepackage{tikz}\usetikzlibrary{cd,fit,matrix,arrows,decorations.pathmorphing}
\tikzset{commutative diagrams/.cd}

\usepackage[sort&compress]{natbib}

\usepackage{cleveref}
\usepackage{xcolor}

\numberwithin{equation}{section}

\newtheorem{theorem}{Theorem}[section]
\newtheorem{corollary}[theorem]{Corollary}
\newtheorem{lemma}[theorem]{Lemma}
\newtheorem{proposition}[theorem]{Proposition}
\newtheorem{conjecture}[theorem]{Conjecture}

\theoremstyle{definition}
\newtheorem{definition}[theorem]{Definition}
\newtheorem{definition-theorem}[theorem]{Definition-Theorem}
\newtheorem{definition-lemma}[theorem]{Definition-Lemma}
\newtheorem{example}[theorem]{Example}

\newtheorem{question}[theorem]{Question}

\newtheorem{remark}[theorem]{Remark}

\theoremstyle{remark}

\newtheorem*{remark*}{Remark}

\usepackage{enumitem}
\setenumerate[1]{label=\textup{(\alph*)}}
\setenumerate[2]{label=\textup{(\roman*)}}

\newcommand\Z{{\mathbb Z}}
\newcommand\N{{\mathbb N}}
\newcommand\Q{{\mathbb Q}}
\newcommand\C{{\mathbb C}}

\def\O{\mathcal O}

\DeclareMathOperator{\im}{im}
\DeclareMathOperator{\coker}{coker}

\DeclareMathOperator{\Hom}{Hom}
\DeclareMathOperator{\End}{End}

\DeclareMathOperator{\Spec}{Spec}

\DeclareMathOperator{\GL}{GL}

\DeclareMathOperator{\Mat}{Mat}

\newcommand\subeq{\subseteq}

\newcommand\supeq{\supseteq}

\newcommand\bbar{\overline} 
\newcommand\tl{\widetilde}
\newcommand\hhat{\widehat} 
\newcommand\onto{\twoheadrightarrow}

\newcommand{\map}[1][]{{\xrightarrow{#1}}} 

\DeclarePairedDelimiter{\set}{\{}{\}}

\DeclarePairedDelimiter{\parens}{\lparen}{\rparen}

\newcommand{\ev}[2]{\left. {#1} \right\rvert_{#2}}

\newcommand{\Zhat}{\hhat{Z}}

\newcommand{\Coh}{\mathrm{Coh}}
\DeclareMathOperator{\Quot}{Quot}
\DeclareMathOperator{\rk}{rk}
\newcommand{\m}{\mathfrak{m}}
\DeclareMathOperator{\Span}{span}
\renewcommand{\L}{\mathbb{L}}

\DeclareMathOperator{\Fil}{Fil}

\DeclareMathOperator{\Gr}{Gr}

\newcommand{\sbp}[1]{_{(#1)}}

\newcommand{\KVar}[1]{K_0(\mathrm{Var}_{#1})}

\newcommand{\KStck}[1]{{K_0(\mathrm{Stck}_{#1})}}

\usepackage{accents}
\newcommand{\ut}{\undertilde}

\newcommand{\calE}{{\mathcal{E}}}
\newcommand{\calF}{{\mathcal{F}}}

\newcommand{\qbinom}[2]{\genfrac{[}{]}{0pt}{}{#1}{#2}}

\begin{document}

\title{Flag Quot and Coh zeta functions}


\author{Yifeng Huang}
\address{Department of Mathematics, Emory University, Atlanta, GA 30322}
\curraddr{Department of Mathematics, Emory University, Atlanta, GA 30322}
\email{yifeng.huang@emory.edu}
\thanks{}

\author{Ruofan Jiang}
\address{Department of Mathematics, University of California--Berkeley}
\curraddr{Department of Mathematics, University of California--Berkeley}
\email{ruofanjiang@berkeley.edu}
\thanks{}


\date{\today}

\dedicatory{}

\begin{abstract}
We study flag Quot and Coh zeta functions on schemes. We first establish a motivic stability result realizing flag Coh zeta functions as infinite-rank limits of flag Quot zeta functions. We then compute the zeta functions for the formal cusp
$R=k[\![T^2,T^3]\!]$. These zeta functions reflect structural properties of the Hall algebra of coherent modules on $R$. In particular, we study a family of Hecke algebras over $R$, and show that their infinite-rank stabilization is commutative.
\end{abstract}

\maketitle

\section{Introduction}
\subsection{Background and motivation}Hilbert schemes of points are important objects that lie at the interface of algebraic geometry, representation theory, topology and physics. In \cite{huangjiang2023torsionfree}, the authors have studied Quot schemes of points and stacks of finite coherent sheaves (which are higher-rank and infinite-rank generalizations of Hilbert schemes of points) from an enumerative perspective. More precisely, we studied the generating function (in the Grothendieck ring of motives) of the  Quot schemes of $n$-points of a torsion-free sheaf, with varying $n$ (called the \textit{Quot zeta function}); as well as  the generating function of stacks of torsion coherent sheaves of length $n$, with varying $n$ (called the \textit{Coh zeta function}). In nice cases (e.g. when the base scheme is a reduced curve), these zeta functions have good properties; cf. \cite{huangjiang2023torsionfree}.

In this paper, we consider the flag generalizations of Quot schemes of points and finite coherent stacks.
The flag objects capture essential information not seen in the single-step version. For example, flag objects parametrize successive ``Hecke modifications'', which are important from a representation-theoretic point of view: a typical relation between flag and non-flag objects is the following
\begin{center}
\begin{tikzcd}
          & {\Coh_{n,n+1}(X)} \arrow[ld, "\stackrel{\leftarrow}{h}"'] \arrow[rd, "\stackrel{\rightarrow}{h}"] &               \\
\Coh_n(X) &                                                                                                   & \Coh_{n+1}(X)
\end{tikzcd}
\end{center}
Here $\Coh_n(X)$ and $\Coh_{n+1}(X)$ are stacks of coherent sheaves of length $n$ and $n+1$, respectively; and $\Coh_{n,n+1}(X)$ is the stack parametrizing a two-step flag. The structure allows the translation of cohomological data between $\Coh_n(X)$ and $\Coh_{n+1}(X)$ through a push--pull process. On the Quot side, there is a similar story. Such structures appear in various works and are closely related to quantum group representations; to name a few, see \cite{Nakajima1,nakajima1999hilbert,Negut1,Negut2,CoHa1}. When $X$ is a smooth curve, the Grothendieck--Springer resolution of $\Coh_n(X)$ is the flag coherent stack $\Coh_{1,2,\dots,n}(X)$ parameterizing full flags; cf. \cite{Laumon87,FYZ1}. In fact, flag objects are analogues of or directly related to shtukas. That said, it would be remiss not to mention that the fascination of flag Quot schemes and coherent stacks extends beyond these connections, in that they can be formulated for arbitrary schemes, which enables one to envision a singular and higher-dimensional analogue of geometric Langlands, as discussed in various works on (cohomological/motivic) Hall algebras; cf.  \cite{kapranov2000,KapEisenstein,CoHa1}.

\subsection{Main results} We will mainly focus on the enumerative aspect of flag Quot schemes and stacks of flag coherent modules, and see how they arise as a shadow of the Hecke philosophy discussed above.

Let us briefly introduce the main objects of our study:
Let $X$ be a formal scheme that is of formal finite type\footnote{This means that $X$ is Noetherian and has a finite-type scheme of definition. In particular, $X$ has finite-type truncations.} over an arbitrary field $k$. Throughout, $\N=\{0,1,2,\ldots\}$. Let $\mathbf{d}=(d_1,\ldots,d_r)$ satisfy $0\leq d_1\leq\cdots\leq d_r$. Let $\mathcal{E}$ be a coherent sheaf on $X$. Then the \textit{flag Quot scheme}  $\mathrm{Quot}_{\mathbf{d}}^R(\mathcal{E})$ classifies descending flags of submodules of $\mathcal{E}\supeq \mathcal{K}_1\supeq\dots \supeq \mathcal{K}_r$ such that $\dim_k\mathcal{E}/\mathcal{K}_i=d_i$.\footnote{Here and hereafter, the notation $\dim_k \mathcal{F}=d$ means $\mathcal{F}$ is a finite-length coherent sheaf with $\dim_k H^0(\mathcal{F})=d$.} The \textit{flag finite coherent stack}  $\mathrm{Coh}_{\mathbf{d}}(X)$ classifies ascending flags of finite coherent sheaves $\mathcal{F}_1\subeq\dots\subeq \mathcal{F}_r$  such that $\dim_k\mathcal{F}_i=d_i$. Equivalently, by taking $\mathcal{Q}_i=\mathcal{E}/\mathcal{K}_i$ or $\mathcal{Q}_i=\mathcal{F}_r/\mathcal{F}_{r-i}$ as appropriate, with $\mathcal{F}_0=0$, $\mathrm{Quot}_{\mathbf{d}}^R(\mathcal{E})$ classifies surjective flags $\mathcal{E}\onto \mathcal{Q}_r\onto\dots\onto \mathcal{Q}_1$ with $\dim_k \mathcal{Q}_i=d_i$, and $\mathrm{Coh}_{\mathbf{d}}(X)$ classifies surjective flags $\mathcal{Q}_r\onto\dots\onto \mathcal{Q}_1$ with $\dim_k \mathcal{Q}_i=d_r-d_{r-i}$, where $d_0=0$.

For $\mathbf{n}=(n_1,\dots,n_r)\in\N^r$, let $\mathbf{d}(\mathbf{n})=(n_1,n_1+n_2,\dots, n_1+n_2+\dots +n_r)$. Define the $r$-step Quot zeta function and the $r$-step Coh zeta function by
\begin{align*}
{Z}^{X}_{\mathcal{E}}(t_1,\dots,t_r)&:=\sum_{\mathbf{n}\in \N^r} [\Quot^X_{\mathbf{d(n)}}(\mathcal{E})]t_1^{n_1}\cdots t_r^{n_r}\in K_0(\mathrm{Var}_k)[\![t_1,\dots,t_r]\!],\\
\Zhat_X(t_1,\dots,t_r)&:=\sum_{\mathbf{n}\in \N^r} [\mathrm{Coh}_{\mathbf{d(n)}}(X)]t_1^{n_1}\cdots t_r^{n_r}\in K_0(\mathrm{Stck}_k)[\![t_1,\dots,t_r]\!].
\end{align*}
For $r=0$, both zeta functions are defined to be $1$. Here $K_0(\mathrm{Var}_k),K_0(\mathrm{Stck}_k)$ are the Grothendieck rings of varieties/stacks over $k$; the latter is in the sense of Definition~\ref{def:k0stck}. When $r=1$, they are the classical objects studied in \cite{huangjiang2023torsionfree}. Some fundamental one-step results in \cite{huangjiang2023torsionfree} extend naturally to the flag setting. First, we have the following motivic stability result:
\begin{theorem}[The flag Quot to Coh formula]\label{thm:nestQC}
  We have
$$\Zhat_X(t_1,\dots,t_r)=\lim_{d\rightarrow\infty} {Z}_{\mathcal{O}_X^d}^{X}(\L^{-d}t_r,\cdots,\L^{-d}t_1).$$
\end{theorem}

For smooth curves, we have the following well-known formulas for these zeta functions. The first part is the main result of \cite{monavariricolfi} after a change of variables, and the second part follows from the first by Theorem~\ref{thm:nestQC}.
\begin{theorem}\label{thm:flag-cl}
    Let $r\geq 0$, let $X$ be a smooth curve, and let $\mathcal{E}$ be a vector bundle of rank $d$ over $X$. Then
 \begin{equation}
 {Z}^{X}_{\mathcal{E}}(t_1,\dots,t_r)=\prod_{j=1}^r \prod_{i=0}^{d-1} Z_X(\L^{i}t_j),\qquad
\Zhat_X(t_1,\dots,t_r)=\prod_{j=1}^r \prod_{i=1}^\infty Z_X(\L^{-i}t_j),
    \end{equation}
    where $Z_X$ is the motivic zeta function of $X$ as in \textup{\cite{kapranov2000}}.
\end{theorem}

The most mysterious part of the story, which lies at the boundary of our knowledge, is when $X$ is singular or higher dimensional. In the following, we will consider the case of (singular) reduced curves. Since $Z^X_\mathcal{E}$ and $\widehat{Z}_X$ satisfy the cut-and-paste property (Proposition~\ref{prop:cut-and-paste}),
it does not hurt to assume $X$ is a singular curve germ. To this end, we will consider the simplest singularity $X=\mathrm{Spf}R$ where
 $$R=k[\![T^2,T^3]\!].$$ Let $\widetilde{R}$ be the normalization of $R$, viewed as an $R$-module. Let $\overline{\mathrm{Bun}}_{\GL_d,R}$ be the stack of torsion-free modules over $R$. It is known that $\overline{\mathrm{Bun}}_{\GL_d,R}$ has exactly $d+1$ points  $\mathcal{E}_a^{(d)}:=\widetilde{R}^{a}\oplus R^{d-a}$ for $0\leq a\leq d$, and the topology is given by the chain of specializations $\mathcal{E}_0^{(d)}\rightarrow \cdots \rightarrow\mathcal{E}_d^{(d)}$.

\begin{definition}\label{def:defcore}\begin{enumerate}\item Let $\Quot_{n;\mathcal{E}_b^{(d)}}^R(\mathcal{E}_a^{(d)})\subseteq \Quot_{n}^R(\mathcal{E}_a^{(d)})$ be the (reduced) locally closed subvariety that pointwise parametrizes codimension $n$ submodules of $\mathcal{E}_a^{(d)}$ that are furthermore isomorphic to $\mathcal{E}_b^{(d)}$.
    \item Let $H^{(d)}_{R,n}\in \Mat_{d+1}(\KVar{k})$ be the $(d+1)\times (d+1)$-matrix with $(a,b)$-th entry $[\Quot_{n;\mathcal{E}_b^{(d)}}^R(\mathcal{E}_a^{(d)})]$ (we count from $0$, the same for below). We call it the \textit{$n$-th Hecke transition matrix}.
\item Let  $\widehat{H}_{R,n}=\lim_{d\rightarrow\infty}\L^{-nd}H^{(d)}_{R,n}$ be the infinite rank matrix whose $(a,b)$-th entry is $$\lim_{d\rightarrow\infty}\L^{-nd}[\Quot_{n;\mathcal{E}_b^{(d)}}^R(\mathcal{E}_a^{(d)})].$$ We call it the \textit{$n$-th stable Hecke transition matrix}.
\item Let $H^{(d)}_{R}(t)=\sum_{n\geq 0}H^{(d)}_{R,n}t^n$ and $\widehat{H}_{R}(t)=\sum_{n\geq 0}\widehat{H}_{R,n}t^n$, so that $\widehat{H}_R(t)=\lim_{d\to\infty} H^{(d)}_R(\L^{-d}t)$ coefficient-wise.
\end{enumerate}
\end{definition}
The constructions given above are closely related to Hall/Hecke  algebras, and will be explained in \S\ref{subsec:hhalgebras}. We use the standard notation
\[ (a;q)_n=(1-a)(1-aq)\cdots (1-aq^{n-1}), \quad 0\leq n\leq \infty\]
and
\[ \qbinom{n}{r}_q=\begin{cases}
    \frac{(q;q)_n}{(q;q)_r(q;q)_{n-r}}\in \N[q], & 0\leq r\leq n;\\
    0, & \text{otherwise}.
\end{cases}\]
\begin{theorem}\label{eq:nest cusp}
    Notation as above.  \begin{enumerate}
        \item The $(a,b)$-th entry of $H_{R}^{(d)}(t)$ is given by  \begin{equation}\label{eq:nest cusp 2}
\sum_{l=\max\{a,b\}}^{d}{\qbinom{d-a}{d-l}_{\L}\qbinom{l}{b}_{\L}}\frac{\L^{(l-b)d}t^{2l-a-b}}{(\L^{d-l}t\,;\,\L)_{l}}.
\end{equation}
\item
The $(a,b)$-th entry of $\widehat{H}_R(t)$ is given by  \begin{equation}\label{eq:limit QUOT}
 \sum_{l\geq \max\{a,b\}}\qbinom{l}{b}_{\mathbb{q}}\frac{\mathbb{q}^{l(l-a)-b(l-b)}t^{2l-a-b}}{(\mathbb{q};\mathbb{q})_{l-a}(\mathbb{q}t;\mathbb{q})_l},\;\;\;\mathbb{q}=\L^{-1}.
\end{equation}
\item ${Z}^{R}_{\mathcal{E}^{(d)}_a}(t_1,\dots,t_r)=e_a^T H^{(d)}_R(t_1)H^{(d)}_R(t_2)\cdots {H}_R^{(d)}(t_r) \,\mathbf{1}$, the sum of entries on the $a$-th row of $H^{(d)}_R(t_1)H^{(d)}_R(t_2)\cdots {H}_R^{(d)}(t_r)$.
\item As the $a=0$ special case, $\widehat{Z}_{R}(t_1,\dots,t_r)=e_0^T\widehat H_R(t_1)\widehat H_R(t_2)\cdots\widehat H_R(t_r)\,\mathbf{1}$.
\item $\widehat{H}_R(t)$ admits the following infinite product expansion
$$\widehat{H}_R(t)=  (\mathbb{q}t;\mathbb{q})^{-1}_{\infty}\prod_{s=0}^{\infty}
 \left((1-\mathbb{q}^{s+1}t+\mathbb{q}^{2s+1}t^2) I+(1-\mathbb{q})\mathbb{q}^st\,\widehat{H}_{R,1}\right).$$
Here $\widehat{H}_{R,1}$ is given explicitly in \eqref{eq:stable-hecke-one}, and convergence of the product is part of the claim.
\item $\widehat{Z}_{R}(t_1,\dots,t_r)$ is symmetric in $t_1,\dots, t_r$.
\end{enumerate}
\end{theorem}
Among these results, parts (e)(f) are the deepest and will be explained further in \S\ref{subsec:hhalgebras}. On the other hand, the following polynomiality result follows readily from (a)(c):
\begin{corollary}
    For $R=k[\![T^2,T^3]\!]$, the normalized flag Quot zeta function $$N^R_{\mathcal{E}}(t_1,\dots,t_r)=\frac{{Z}^{R}_{\mathcal{E}}(t_1,\dots,t_r)}{{Z}^{\tl R}_{{\tl R}^d}(t_1,\dots,t_r)}$$ is a polynomial in $\KVar{k}[t_1,\dots,t_r]$ for any torsion-free $R$-module $\mathcal{E}$ of rank $d$.
\end{corollary}
This is a case of the following \begin{conjecture}
   For any reduced singular curve germ $R$, the normalized flag Quot zeta function
$N^R_{\mathcal{E}}(t_1,\dots,t_r)$ as defined above is a polynomial in $\KVar{k}[t_1,\dots,t_r]$.
\end{conjecture}  When $r=1$, this is proved in \cite{huangjiang2023torsionfree}, and the subsequent works \cite{HJO,hlop2026,lauono2026,hlo2026} demonstrate that its computation is the key to unlocking many surprising connections. Strong evidence for the conjecture is that when $k=\mathbb{F}_q$, $N^R_{\mathcal{E}}(t_1,\dots,t_r)$ is indeed a polynomial in $\Z[t_1,\dots,t_r]$ upon passing to the point count measure. This is an exercise based on the rationality result on $p$-adic integration \cite{bushnellreiner1980zeta,knopfmacher1985}.

Another question, to which we don't have enough evidence, is \begin{question}
     For a reduced singular curve germ $R$, when is $\widehat{Z}_{R}(t_1,\dots,t_r)$ symmetric?
\end{question}

\begin{remark}
    By Theorem~\ref{eq:nest cusp}, the motive of the full flag Coh stack $\Coh_{(1,\dots,d)}(R)$ can be determined in principle. To provide an independent consistency check, we present a different method to compute it by classifying commuting strictly upper triangular matrices $(A,B)$ such that $AB=BA, A^2=B^3$, which is interesting on its own. The details are provided in the Appendix. The basis of the argument appears in our unpublished note \cite{huangjiang2023punctual}; to make the paper self-contained, we reproduce the proof in full.
\end{remark}
\subsection{Hall and Hecke algebras}\label{subsec:hhalgebras} Let $k=\mathbb{F}_q$ be a finite field and let $\mathrm{Hall}(R)$
be the Hall algebra of coherent $R$-modules\footnote{We don't need to restrict $k$ if we use motivic Hall algebras.}. Let's recall the definition: let $|\mathrm{Coh}(R)|$ be the set of isomorphism classes of coherent $R$-modules. Then $\mathrm{Hall}(R)$ is the $\C$-vector space of finitely supported functions $|\mathrm{Coh}(R)|\rightarrow \C$ with product structure
$$(f\circ g)(A) = \sum_{A'\subseteq A} f(A')g(A/A').$$
In $\mathrm{Coh}(R)$ there are two families of substacks: the stack of dimension $n$ finite coherent modules $\mathrm{Coh}_n(R)$, and the stack of rank $d$ torsion-free modules $\overline{\mathrm{Bun}}_{\GL_d,R}$. Let $\mathrm{Hall}_{\mathbb{N}}(R)$ be the sub-algebra consisting of functions supported on the substack of finite coherent modules $\mathrm{Coh}_{\mathbb{N}}(R)=\bigsqcup_{n\geq 0}\mathrm{Coh}_n(R)$.

For each $d\geq 0$, let $\mathcal{A}_{\GL_d,R}$ be the space of $\mathbb{C}$-valued functions on the set $|\overline{\mathrm{Bun}}_{\GL_d,R}|$. This is a local and singular analogue of automorphic forms. Observe that we have a right action of $\mathrm{Hall}_{\mathbb{N}}(R)$ on $\mathcal{A}_{\GL_d,R}$ derived from the product structure on $\mathrm{Hall}(R)$. It can be described as follows: let $\mathbf{1}_A$ be the characteristic function for an $A\in |\Coh(R)|$. Then for $A\in |\Coh_{\mathbb{N}}(R)|$, the action is
$$ \mathbf{1}_{\mathcal{E}^{(d)}_a} \circ \mathbf{1}_A = \sum_{b=0}^{d}\#\left\{K\subseteq \mathcal{E}^{(d)}_a|K\simeq \mathcal{E}^{(d)}_b,\mathcal{E}^{(d)}_a/K\simeq A\right\} \cdot \mathbf{1}_{\mathcal{E}^{(d)}_b}.$$
Let $\mathbf{1}_{n}$ be the sum of characteristic functions for elements in $|\Coh_{n}(R)|$. Then the action of $\mathbf{1}_n$ on $\mathcal{A}_{\GL_d,R}$ is an analogue of a Hecke action by ``modifying $n$ points down''. We will denote this action by $H_{R,n}^{(d)}$ and call it the \textit{$n$-th Hecke operator}. With respect to the basis  $\{\mathbf{1}_{\mathcal{E}^{(d)}_a}\}_{a=0,1,\dots,d}$ of $\mathcal{A}_{\GL_d,R}$, $H_{R,n}^{(d)}$ is precisely presented by the $n$-th Hecke transition matrix in Definition~\ref{def:defcore}.
\begin{definition}
Define the (possibly non-commutative) algebra $\mathcal{H}^{(d)+}_R$ (resp. $\widehat{\mathcal{H}}_R^+$) to be the algebra generated over $\Q(q)$\footnote{If the reader takes a point-count perspective, then $\Q(q)=\Q$ as $q\in \mathbb{N}$. If the reader is more motivically minded, then they can treat $q$ as a formal variable (playing the same role as $\L$).}  by $H^{(d)}_{R,1},H^{(d)}_{R,2},\dots$ (resp. $\widehat{H}_{R,1},\widehat{H}_{R,2},\dots $). It is called the (positive) \textit{rank $d$ Hecke algebra} (resp. \textit{stable Hecke algebra}) for $R$.
\end{definition}
Note that $\widehat{\mathcal{H}}_R^+$ acts on the vector space $\varprojlim_d \mathcal{A}_{\GL_d,R}=\prod_{a\geq 0}\mathbb{C}\cdot\mathbf{1}_{\mathcal{E}^{(\infty)}_{a}}$ (limit taken with transition morphisms induced by $\mathcal{E}^{(d)}_a\rightarrow \mathcal{E}^{(d+1)}_a$). In general, $\mathcal{H}^{(d)+}_R$ is non-commutative, but surprisingly
\begin{theorem}\label{thm:commuteHecke} $\widehat{\mathcal{H}}_{R}^+=\Q(q)[\widehat{H}_{R,1}]$, hence commutative.

\end{theorem}

In fact, this is a consequence of the following three-term recursion (Theorem~\ref{thm:recursion}):
\begin{equation}\label{eq:recursion_2in q}
\widehat{H}_{R,1}\cdot\widehat{H}_{R,m}
 =[m+1]_q\widehat{H}_{R,m+1}-[m]_q\widehat{H}_{R,m}
 -\frac q{q-1}\widehat{H}_{R,m-1}.
\end{equation} While Theorem~\ref{thm:commuteHecke} is the structural reason behind Theorem~\ref{eq:nest cusp}(e)(f), one can further question the structural reason behind the formula (\ref{eq:recursion_2in q}). It turns out that this is closely related to a certain quantum algebra action, which we will explore in a subsequent work.

\subsection{One-step Quot zeta function}
\label{subsec:intro-one-step-interpolation}

We note that the one-step consequence of Theorem~\ref{eq:nest cusp} is
already interesting.  In \cite{huangjiang2023torsionfree}, the authors
computed the normalized Quot zeta functions associated with the two extreme
modules $R^d$ and $\widetilde R^d$, and curiously noted that both have simple
finitized, $t$-deformed Rogers--Ramanujan expressions that differ only by the
substitution $t\mapsto t^2$.  The intermediate torsion-free modules
\[
  \mathcal E_a^{(d)}=\widetilde R^{a}\oplus R^{d-a},
  \qquad 0\leq a\leq d,
\]
provide a natural interpolation between these two formulas.

Recall that
\[
 N^R_{\mathcal E_a^{(d)}}(t)
 :=\frac{Z^R_{\mathcal E_a^{(d)}}(t)}
 {Z^{\widetilde R}_{\widetilde R^d}(t)}
 =(t;\L)_d Z^R_{\mathcal E_a^{(d)}}(t).
\]
The one-step specialization of the main theorem admits the following
positive formula; see Section~\ref{sec:one-step-interpolation} for the proof.

\begin{corollary}\label{cor:intro-one-step-interpolation}
For $0\leq a\leq d$, one has
\begin{equation}\label{eq:intro-one-step-interpolation}
 N^R_{\mathcal E_a^{(d)}}(t)
 =
 \sum_{i=0}^{a}\sum_{j=0}^{d-a}
 \L^{d(i+j)+(a-i)j}
 \qbinom{a}{i}_{\L}
 \qbinom{d-a}{j}_{\L}
 t^{i+2j}.
\end{equation}
In particular,
\begin{align*}
 N^R_{R^d}(t)
 =\sum_{j=0}^{d}\qbinom{d}{j}_{\L}\L^{dj}t^{2j},\qquad
 N^R_{\widetilde R^d}(t)
 =\sum_{i=0}^{d}\qbinom{d}{i}_{\L}\L^{di}t^i,
\end{align*}
and hence
\[
 N^R_{R^d}(t)=N^R_{\widetilde R^d}(t^2).
\]
\end{corollary}

There is also a normalization of
\eqref{eq:intro-one-step-interpolation} that admits a stable limit.  Put $\mathbb q=\L^{-1}$ and define
\[
 \mathscr R_{r,s}(x,y;\mathbb q)
 :=\sum_{i=0}^{r}\sum_{j=0}^{s}
 \mathbb q^{i^2+ij+j^2}
 \qbinom{r}{i}_{\mathbb q}
 \qbinom{s}{j}_{\mathbb q}x^iy^j;
\]
note that one or both of $r,s$ can be $\infty$, as $$\qbinom\infty i_q\coloneqq \lim_{n\to \infty}\qbinom ni_q=\frac{1}{(q;q)_\infty}.$$
Then
\begin{equation}\label{eq:intro-a2-slice}
 N^R_{\mathcal E_a^{(d)}}(\L^{-d}t)
 =\mathscr R_{a,d-a}(\mathbb q^{-a}t,t^2;\mathbb q),
\end{equation}
so that
\begin{equation}
 \lim_{d\to \infty}N^R_{\mathcal E_a^{(d)}}(\L^{-d}t)
 =\mathscr R_{a,\infty}(\mathbb q^{-a}t,t^2;\mathbb q),
\end{equation}
a two-fold infinite sum. At $a=0$, \eqref{eq:intro-a2-slice} reduces to the one-fold sum that appears in \cite{huangjiang2023torsionfree}
\[
 N^R_{R^d}(\L^{-d}t)
 =\sum_{j=0}^{d}\mathbb q^{j^2}
   \qbinom{d}{j}_{\mathbb q}t^{2j},
\]
whose $t=1$ specialization in the limit $d\to\infty$ is one of the two Rogers--Ramanujan identities:
\[ \lim_{d\to \infty} N^R_{R^d}(\L^{-d})=\sum_{j=0}^\infty \frac{\mathbb q^{j^2}}{(\mathbb q;\mathbb q)_j}=\frac{1}{(\mathbb q;\mathbb q^5)_\infty(\mathbb q^4;\mathbb q^5)_\infty}.\]

\subsection*{AI disclosure} We used ChatGPT-5.6 Sol to discover the key recursion (\ref{eq:recursion_2in q}) and a positive formula \eqref{eq:intro-one-step-interpolation}, and we are responsible for their correctness. Other computations and the interpretation of (\ref{eq:recursion_2in q}) in terms of Hall/Hecke algebras belong to the authors.

\section{Main players}
We review the definition of classical Quot schemes of points, finite coherent stacks and commuting varieties. We then make flag generalizations.
    \subsection{Classical objects} Let $X$ be a quasi-projective scheme over $k$. Let $d\geq 0$ be an integer.
\begin{definition}[Quot schemes of points]  Given a coherent sheaf $\mathcal{E}$ on $X$, define
the Quot scheme of $d$-points $\mathrm{Quot}_{d}^R(\mathcal{E})$ to be
the functor sending an affine $k$-scheme $S$ to $\mathcal{O}_{X_S}$-submodules of $\pi^* \calE$ whose quotients are locally free of rank $d$ over $S$; here $X_S:=S\times_k X$ and $\pi:X_S\to X$ is the structural morphism. It is classically known to be represented by a  scheme over $k$.
\end{definition}
\begin{definition}[Finite coherent stacks] Define $\Coh_d(X)$ to be the stack of finite coherent sheaves over $X$. As a functor, it sends an affine $k$-scheme $S$ to the groupoid of finite $\mathcal{O}_{X_S}$-modules that are flat of rank $d$ over $S$. It is classically known to be an Artin stack over $k$.
\end{definition}
\begin{definition}[Commuting varieties]
   Define the commuting variety to be the  following functor: for a $k$-algebra $S$,
$$\mathrm{C}_{d}(X)(S)=\left\{ (\mathcal{F},\iota)\ \vrule \ \begin{array}{l} \,  \mathcal{F}\in \mathrm{Coh}_d(X)(S), \,\iota: k^d\otimes S\xrightarrow{\sim} H^0(\mathcal{F})\end{array}\right\}.$$
\end{definition}
There is a natural $\GL_d$ action on $\mathrm{C}_{d}(X)$ permuting the framing $\iota$. The functor is representable because of the following matrix point interpretation in the affine case:
\begin{proposition}[Matrix points]
Present $R$ by $R=k[x_1,\cdots,x_n]/(f_1,\cdots,f_m)$. Then $\mathrm{C}_{d}(R)$ is the scheme classifying ``matrix points of $R$'' in the following sense:
\begin{equation}\mathrm{C}_{d}(R):=
\set*{
X_1,\cdots,X_n\in \Mat_d
\ \vrule \
\begin{array}{l}
\forall j, f_j(X_1,\cdots, X_n)=0, \\\forall i,i', [X_i,X_{i'}]=0.
\end{array}
}.
\end{equation}
The definition does not depend on the presentation.
\end{proposition}

\begin{corollary}
We have
$ \mathrm{Coh}_{d}(X)=[\mathrm{C}_{d}(X)/\GL_d]$. Since étale $\GL_d$-torsors are Zariski, we have $[\mathrm{Coh}_{d}(X)]\cdot [\GL_d]=[\mathrm{C}_{d}(X)]$ in $K_0(\mathrm{Stck}_k)$.
\end{corollary}
A proof of these facts can be found in \cite{huangjiang2023torsionfree}.
\subsubsection{Ind-generalization}\label{subsubsec:ind} We will also need generalizations of these objects to ind-schemes, in particular for formal schemes. Let $X$ be an ind-finite-type ind-scheme over $k$ presented as a limit $\bigcup_{n\geq 0} X_n$, and let $\mathcal{E}$ be an ind-sheaf on $X$. We define
\begin{align*}
\Quot^X_{d}(\mathcal{E})&=\varinjlim \Quot_{d}^{X_n}(\mathcal{E}|_{X_n} ),\\
\mathrm{Coh}_{d}(X)&=\varinjlim \mathrm{Coh}_{d}(X_n),\\
\mathrm{C}_{d}(X)&=\varinjlim \mathrm{C}_{d}(X_n).
\end{align*}
It is straightforward to check that $\Quot^X_{d}(\mathcal{E})$, $\mathrm{Coh}_{d}(X)$ and $\mathrm{C}_{d}(X)$ are ind-finite-type ind-objects that do not depend on the presentation of $X$. Furthermore, suppose that $X$ is a Noetherian formal scheme with a finite-type scheme of definition. Choose an ideal of definition $\mathcal I$ and take $X_n=(|X|,\mathcal O_X/\mathcal I^{n+1})$. The reduced structures of $\Quot^{X_n}_d(\mathcal E|_{X_n})$, $\mathrm{Coh}_d(X_n)$ and $\mathrm C_d(X_n)$ stabilize for $n\ge d$. Indeed, for a length-$d$ sheaf $\mathcal F$ representing a geometric point of any truncation, the filtration $\mathcal F\supseteq\mathcal I\mathcal F\supseteq\cdots$ strictly decreases until zero by Nakayama's lemma, so $\mathcal I^d\mathcal F=0$. The transition maps between these moduli objects are closed immersions and, from this stage onward, are surjective on geometric points; hence they induce isomorphisms on reductions. In particular, the reduced ind-objects above are finite-type schemes or stacks; cf.~\cite[\S 2]{huangjiang2023torsionfree} and \cite[\S 2.2.1]{fantechiricolfi2024motivic}.
\subsection{Flag objects}\label{subsec:flag-objects}Let $X$ be a quasi-projective scheme over $k$. Let $\mathbf{d}=(d_1,\ldots,d_r)$ be a weakly increasing sequence of nonnegative integers. This will be called a \textit{dimension vector}. Fix a flag of $k$-vector spaces
 $\mathbf{V}:=(V_1\subseteq \cdots\subseteq V_r)$ such that $\dim V_i=d_i$.  Let $\mathrm{P}(\mathbf{V})\subseteq \GL(V_r)$ be the parabolic fixing $\mathbf{V}$.
We use $V_0=0$; the empty flag has moduli space $\Spec k$.
\begin{definition}[Flag Quot schemes of points]\label{def:def Nest Quot} Define
the \textit{$\mathbf{d}$-flag Quot scheme of points} $\mathrm{Quot}_{\mathbf{d}}^R(\mathcal{E})$ of a coherent sheaf $\mathcal{E}$ over $X$ to be the closed subscheme of $\prod_{i=1}^r\mathrm{Quot}_{d_i}^R(\mathcal{E})$ defined by the incidence relation $\mathcal{K}_1\supseteq \mathcal{K}_2\supseteq \cdots \supseteq \mathcal{K}_{r}$,
 where $\mathcal{K}_i$ is the submodule determined by a point in  $\mathrm{Quot}_{d_i}^R(\mathcal{E})$.
 \end{definition}

\begin{definition}[Flag finite coherent stacks]
The \textit{$\mathbf d$-flag finite coherent stack} $\mathrm{Coh}_{\mathbf d}(X)$ sends an affine $k$-scheme $S$ to the groupoid of flags
\[
0=\mathcal F_0\subseteq\mathcal F_1\subseteq\cdots\subseteq\mathcal F_r
\]
of coherent sheaves on $X_S$ with support finite over $S$, such that each $\mathcal F_i/\mathcal F_{i-1}$ is $S$-flat of relative length $d_i-d_{i-1}$, where $d_0=0$. Morphisms are isomorphisms of flags.
\end{definition}

\begin{definition}[Staircase commuting varieties] Define
the \textit{$\mathbf{d}$-staircase commuting variety} to be the  following functor: for a $k$-algebra $S$,
$$\mathrm{C}_{\mathbf{d}}(X)(S)=\left\{ (\mathcal{F}_1,\cdots, \mathcal{F}_r,\iota)\ \vrule \ \begin{array}{l} \,  \mathcal{F}_i\in \mathrm{Coh}_{d_i}(X)(S), \,\iota: V_{r}\otimes S\xrightarrow{\sim} H^0(\mathcal{F}_{r})\\ \text{such that }\iota |_{V_i\otimes S}:  V_{i}\otimes S\xrightarrow{\sim} H^0(\mathcal{F}_{i}).\end{array}\right\}.$$
\end{definition}
There is a natural transformation $\mathrm{C}_{\mathbf{d}}(X)\rightarrow \mathrm{C}_{d_i}(X)$ for each $i$ which, over a $k$-algebra $S$, is the map $(\mathcal{F}_1,\cdots, \mathcal{F}_r,\iota)\rightarrow (\calF_i, \iota|_{V_i\otimes S})$.
There is a natural action of $\mathrm{P}(\mathbf{V})$ on $\mathrm{C}_{\mathbf{d}}(X)$, changing the framing $\iota$.

It is also not too hard to see the following matrix point interpretation in the affine case:
\begin{proposition}[Parabolic matrix points]\label{prop:parabolicma}
Let $R=k[x_1,\cdots,x_n]/(f_1,\cdots,f_m)$. Then $\mathrm{C}_{\mathbf{d}}(R)$ is the scheme classifying ``parabolic matrix points of $R$'' in the following sense:
\begin{equation}\label{eq:definmatrix_flag}\mathrm{C}_{\mathbf{d}}(R):=
\set*{
X_1,\cdots,X_n\in \mathrm{Lie}\,\mathrm{P}({\mathbf{V}})
\ \vrule \
\begin{array}{l}
\forall j, f_j(X_1,\cdots, X_n)=0, \\\forall i,i', [X_i,X_{i'}]=0.
\end{array}
}.
\end{equation}

The definition does not depend on the presentation.
\end{proposition}

\begin{corollary}\label{cor:artin}
We have
$ \mathrm{Coh}_{\mathbf{d}}(X)=[\mathrm{C}_{\mathbf{d}}(X)/\mathrm{P}(\mathbf{V})]$. Since étale $\mathrm{P}(\mathbf{V})$-torsors are Zariski, we have $[\mathrm{Coh}_{\mathbf{d}}(X)]\cdot [\mathrm{P}(\mathbf{V})]=[\mathrm{C}_{\mathbf{d}}(X)]$ in $K_0(\mathrm{Stck}_k)$.
\end{corollary}
The proof of these facts parallels the classical theory.

\medskip

The proof of Theorem~\ref{thm:nestQC} requires comparing the flag objects in the equivalent quotient descriptions, so we record them below. Over an affine base $S$, the scheme $\Quot^X_{\mathbf d}(\mathcal E)$ parametrizes surjective flags
\[
\mathcal E_S\onto\mathcal Q_r\onto\cdots\onto\mathcal Q_1,
\qquad \mathcal Q_i\text{ is $S$-flat of relative length }d_i,
\]
where $\mathcal E_S$ is the pullback of $\mathcal E$. On the Coh side, put $\widehat d_i=d_r-d_{r-i}$ and $\widehat{\mathbf d}=(\widehat d_1,\ldots,\widehat d_r)$. The mutually inverse constructions
\[
\mathcal Q_i=\mathcal F_r/\mathcal F_{r-i},
\qquad
\mathcal F_i=\ker(\mathcal Q_r\onto\mathcal Q_{r-i}),
\qquad \mathcal F_0=\mathcal Q_0=0,
\]
identify $\Coh_{\mathbf d}(X)(S)$ with the groupoid of surjective flags $\mathcal Q_r\onto\cdots\onto\mathcal Q_1$ whose terms are $S$-flat of relative lengths $\widehat d_i$. The framing in $\mathrm C_{\mathbf d}(X)$ becomes a compatible collection of isomorphisms
\[
(V_r/V_{r-i})\otimes_k\mathcal O_S\xrightarrow{\sim}p_*\mathcal Q_i,
\qquad p:X_S\to S.
\]
Such compatible frames form a Zariski-locally trivial $\mathrm P(\mathbf V)$-torsor. In particular, adding a surjection $\mathcal E_S\onto\mathcal Q_r$ to a quotient flag representing $\Coh_{\mathbf d}(X)$ gives precisely an object of $\Quot^X_{\widehat{\mathbf d}}(\mathcal E)$.

\subsubsection{Ind-generalization}
The discussion of \S\ref{subsubsec:ind} applies equally to the flag objects, with stabilization for $n\ge d_r$. Since we will be using it later, we make the following remark
\begin{remark}\label{rmk:indformal}
    If $X=\mathrm{Spf} \,R$ is a formal scheme where $R=k[\![x_1,\cdots,x_n]\!]/(f_1,\cdots,f_m)$, then $\mathrm{C}_{\mathbf{d}}(X)_{\mathrm{red}}=\mathrm{C}_{\mathbf{d}}(X_{N})_{\mathrm{red}}$ for some finite-type truncation $X_{N}\subseteq X$. Therefore, $\mathrm{C}_{\mathbf{d}}(X)_{\mathrm{red}}$ is the reduced structure of the parabolic matrix points (\ref{eq:definmatrix_flag}) with the further assumption that all $X_i$'s are nilpotent.
\end{remark}

\subsubsection{Refinements}\label{subsub:ref1} Let us assume that $X$ is of finite type for technical simplicity. We can refine $\Quot_{\mathbf{d}}^{R}(\calE)$  by requiring $\mathcal{K}_i$ in Definition~\ref{def:def Nest Quot} to lie in prescribed isomorphism classes. Similarly, we can refine $\Coh_{\mathbf{d}}(X)$ by requiring $\mathcal{F}_i$ to lie in prescribed isomorphism classes.

The case for $\Coh_{\mathbf{d}}(X)$ can be made precise relatively easily. Let $\mathcal{C}_i \subseteq \Coh_{d_i}(X)$ be a constructible subset, which plays the role of a prescribed range of isomorphism classes. Let $\mathcal{C}=(\mathcal{C}_1,\dots,\mathcal{C}_r)$. Define the constructible set  $$\Coh_{\mathbf{d};\mathcal{C}}(X):=\Coh_{\mathbf{d}}(X)\cap \bigcap_{i=1}^r\pi_i^{-1}(\mathcal{C}_i),$$ where $\pi_i:\Coh_{\mathbf{d}}(X)\rightarrow \Coh_{d_i}(X)$ is the canonical map. Then $\Coh_{\mathbf{d};\mathcal{C}}(X)$ pointwise classifies flag finite coherent sheaves whose $i$-th coherent sheaf in the flag lies in the prescribed class $\mathcal{C}_i$.

We can also run a similar construction for Quot. Let $\Coh(X)$ be the stack of coherent sheaves. Then we have canonical morphisms $\theta_i:\Quot_{d_i}^{X}(\calE)\rightarrow \Coh(X)$ sending $\mathcal{K}_i$ to its class in $\Coh(X)$. Let $\pi_i: \Quot_{\mathbf{d}}^{X}(\calE)\rightarrow \Quot_{d_i}^{X}(\calE)$ be the morphism sending a flag to $\mathcal{K}_i$. Let $\mathcal{C}_i \subseteq \Coh(X)$ be a constructible subset. Let $\mathcal{C}=(\mathcal{C}_1,\dots,\mathcal{C}_r)$. Define a constructible set  $$\Quot^X_{\mathbf{d};\mathcal{C}}(\mathcal{E}):=\Quot^X_{\mathbf{d}}(\mathcal{E})\cap \bigcap_{i=1}^r\pi_i^{-1}\theta_i^{-1}(\mathcal{C}_i),$$
which is the desired object.
\subsection{Flag zeta functions} Let us restrict to the case where $X$ is a formal scheme of formal finite type over $k$. Let $\mathbf{d},\mathbf{V}$ etc, be the same as in the last subsection. We will always use $\mathbf{n}=\mathbf{n}(\mathbf{d})$ to denote the \textit{grading vector} associated to $\mathbf{d}$, defined as
$$\mathbf{n}=(d_1,d_2-d_1,d_3-d_2,\dots, d_r-d_{r-1}).$$
We will often index by $\mathbf{n}$. We write entries of $\mathbf{n}$ as $(n_1,\cdots, n_r)$, then the associated dimension vector $\mathbf{d}=\mathbf{d}(\mathbf{n})$ is just $(n_1,n_1+n_2,\dots, n_1+n_2+\dots +n_r)$.
\begin{definition}
 Define the \textit{motivic $r$-step flag Quot zeta function} to be the following generating series
    \begin{equation}\label{eq:cl-def-quot}
 {Z}^{X}_{\mathcal{E}}(t_1,\dots,t_r):=\sum_{\mathbf{n}\in \N^r} [\Quot^X_{\mathbf{d}}(\mathcal{E})]t_1^{n_1}\cdots t_r^{n_r}\in K_0(\mathrm{Var}_k)[\![t_1,\dots,t_r]\!].
\end{equation}
\end{definition}
\begin{definition}
 Define the \textit{motivic $r$-step
 flag Coh zeta function} to be the following generating series
    \begin{equation}\label{eq:cl-def-coh}
\Zhat_X(t_1,\dots,t_r):=\sum_{\mathbf{n}\in \N^r} [\mathrm{Coh}_{\mathbf{d}}(X)]t_1^{n_1}\cdots t_r^{n_r}\in K_0(\mathrm{Stck}_k)[\![t_1,\dots,t_r]\!].
\end{equation}
\end{definition}
In future we denote $t^{\mathbf{n}}:=t_1^{n_1}\cdots t_r^{n_r}$ in short.
\begin{proposition}[Cut and paste property]
Let $x\in |X|$ be a closed point, and let $X^{/x}$ be the completion of $X$ at $x$. Then we have factorizations
 \begin{equation}
 {Z}^{X}_{\mathcal{E}}(t_1,\dots,t_r)= {Z}^{X^{/x}}_{\mathcal{E}}(t_1,\dots,t_r)\cdot {Z}^{X\setminus\{x\}}_{\mathcal{E}}(t_1,\dots,t_r)
    \end{equation}
    \begin{equation}
\Zhat_X(t_1,\dots,t_r)=\Zhat_{X^{/x}}(t_1,\dots,t_r)\cdot\Zhat_{X\setminus\{x\}}(t_1,\dots,t_r).
    \end{equation}
    \label{prop:cut-and-paste}
\end{proposition}
The proof is similar to the classical ones and relies on the key fact that finite coherent sheaves on $X$ are supported on finitely many closed points on $X$.
\begin{remark}
    When $k=\mathbb{F}_q$, passing to the point count measure $$\sharp: K_0(\mathrm{Stck}_k)\rightarrow \mathbb{Q},$$ we furthermore have the \textit{Euler product formula}:
     \begin{equation}
\sharp {Z}^{X}_{\mathcal{E}}(t_1,\dots,t_r)=\prod_{x \text{ closed pt in }X}\sharp {Z}^{X^{/x}}_{\mathcal{E}}(t_1,\dots,t_r).
    \end{equation}
    \begin{equation}
\sharp\Zhat_X(t_1,\dots,t_r)=\prod_{x \text{ closed pt in }X}\sharp\Zhat_{X^{/x}}(t_1,\dots,t_r).
    \end{equation}
\end{remark}

\subsubsection{Limits}\label{subsubsec:limit}
The coefficient-wise limits
$$\lim_{r\rightarrow \infty} {Z}^{X}_{\mathcal{E}}(t_1,\dots,t_r),\qquad\lim_{r\rightarrow \infty}\Zhat_X(t_1,\dots,t_r),$$
 converge to elements in $\KVar{k}[\![t_1,t_2,\dots]\!]$ and $\KStck{k}[\![t_1,t_2,\dots]\!]$, respectively. This is due to the elementary observation that deleting redundant steps (where inclusion is equality) does not change the moduli space. We denote the limits by $
 {Z}^{X}_{\mathcal{E}}(t_1,t_2,\dots)$ and  $\Zhat_X(t_1,t_2,\dots)$. These series contain the information of arbitrarily long flags. By the same observation, $
 {Z}^{X}_{\mathcal{E}}(t_1,t_2,\dots)$ and  $\Zhat_X(t_1,t_2,\dots)$ are quasi-symmetric in the variables (in the classical sense of \cite{gesselquasisymmetric} except that they could contain terms of unbounded degrees).
\subsubsection{Refinements}\label{ref:2}
In the situation of \S\ref{subsub:ref1}, we also define zeta functions with prescribed isomorphism classes:
\begin{align*}
{Z}^{X}_{\mathcal{E};\mathcal{C}}(t_1,\dots,t_r)&:=\sum_{\mathbf{n}\in \N^r} [\Quot^X_{\mathbf{d};\mathcal{C}}(\mathcal{E})]t_1^{n_1}\cdots t_r^{n_r},\\
\quad\Zhat_{X,\mathcal{C}}(t_1,\dots,t_r)&:=\sum_{\mathbf{n}\in \N^r} [\mathrm{Coh}_{\mathbf{d};\mathcal{C}}(X)]t_1^{n_1}\cdots t_r^{n_r}.
\end{align*}

\section{The flag Quot to Coh formula}
In this section we prove Theorem~\ref{thm:nestQC}. Following the idea of \cite[\S 3]{huangjiang2023torsionfree}, we view the flag Quot scheme as the open locus of generating sections in a vector bundle over the flag Coh stack. The complement becomes negligible as the number of sections tends to infinity.
We first recall the dimension filtration.

\begin{definition}[\cite{behrenddhillon2007}]
    For $m\in \Z$, let $\Fil_m\KVar{k}[\L^{-1}]$ be the abelian subgroup of $\KVar{k}[\L^{-1}]$ generated by elements of the form $[Y]/\L^n$ where $Y$ is a variety and $\dim Y - n \le m$. This defines a filtration on $\KVar{k}[\L^{-1}]$ satisfying $\Fil_m \cdot \Fil_n\subeq \Fil_{m+n}$, and a sequence $(a_d)_d$ converges to $0$ if $a_d \in \Fil_{m_d}$ with $m_d \to -\infty$ as $d\to\infty$. We note that $[\GL_n]^{-1} \in \Fil_{-n^2}$ and $\L^{-1} \in \Fil_{-1}$. We use $\KStck{k}$ to denote the completion of $\KVar{k}[\L^{-1}]$, with dimension filtration $\Fil_m$ defined as the closure of $\Fil_m\KVar{k}[\L^{-1}]$. \label{def:k0stck}
\end{definition}

Fix $\mathbf d$ and $\mathbf V$ as in \S\ref{subsec:flag-objects}, and consider the parabolic subgroup $P=\mathrm P(\mathbf V)$. All moduli objects in the argument below carry their reduced structures. By the quotient formulation in \S\ref{subsec:flag-objects}, an object of $\Coh_{\mathbf d}(X)$ is a flag
\[
\mathcal Q_r\onto\cdots\onto\mathcal Q_1,
\qquad \dim\mathcal Q_i=\widehat d_i=d_r-d_{r-i}.
\]
For $N\ge0$, let $\mathcal V_{N,\mathbf d}$ parametrize such flags together with a map $\Phi:\mathcal O_X^N\to\mathcal Q_r$. Forgetting $\Phi$ makes $\mathcal V_{N,\mathbf d}$ a vector bundle of rank $Nd_r$ over $\Coh_{\mathbf d}(X)$: in a family over $S$, the map $\Phi$ is an $N$-tuple of sections of the rank-$d_r$ vector bundle $p_*\mathcal Q_r$.
We call $\Phi$ \textit{stable} if it is surjective. The stable locus $\mathcal U_{N,\mathbf d}$ is open, and the quotient description gives a canonical isomorphism
\begin{equation}\label{eq:quot-stable-locus}
\mathcal U_{N,\mathbf d}\simeq\Quot^X_{\widehat{\mathbf d}}(\mathcal O_X^N).
\end{equation}
Pulling back along the compatible-frame bundle $\mathrm C_{\mathbf d}(X)\to\Coh_{\mathbf d}(X)$ trivializes $\mathcal V_{N,\mathbf d}$ as $\mathrm C_{\mathbf d}(X)\times V_r^N$. Over the stable locus, this is precisely the $P$-torsor that forgets the compatible frame of a quotient flag.

\begin{lemma}\label{lem:unstable-locus}
For $\mathcal Z_{N,\mathbf d}:=\mathcal V_{N,\mathbf d}\setminus\mathcal U_{N,\mathbf d}$, we have
\[
[\mathcal Z_{N,\mathbf d}]
\in\Fil_{\dim\mathrm C_{\mathbf d}(X)-\dim P+(N+1)(d_r-1)}.
\]
\end{lemma}
\begin{proof}
Let $\widetilde{\mathcal Z}_{N,\mathbf d}$ be the pullback of $\mathcal Z_{N,\mathbf d}$ to $\mathrm C_{\mathbf d}(X)$, and let $W_N\subseteq V_r^N$ be the determinantal locus of $N$-tuples that do not span $V_r$. Since sections spanning $H^0(X,\mathcal Q_r)$ generate $\mathcal Q_r$ as an $\mathcal O_X$-module, we have
\[
\widetilde{\mathcal Z}_{N,\mathbf d}\subseteq\mathrm C_{\mathbf d}(X)\times W_N.
\]
The bound $\dim W_N\le(N+1)(d_r-1)$ therefore gives
$\dim\widetilde{\mathcal Z}_{N,\mathbf d}\le\dim\mathrm C_{\mathbf d}(X)+(N+1)(d_r-1)$.
Now $[\mathcal Z_{N,\mathbf d}]=[\widetilde{\mathcal Z}_{N,\mathbf d}]/[P]$, and the Levi decomposition of $P$ gives $[P]^{-1}\in\Fil_{-\dim P}$. The result follows.
\end{proof}

\begin{lemma}\label{thm:A-effective}
With the notation above, we have
\[
[\Coh_{\mathbf d}(X)]
-[\Quot^X_{\widehat{\mathbf d}}(\mathcal O_X^N)]\L^{-Nd_r}
\in\Fil_{-N+\dim\mathrm C_{\mathbf d}(X)-\dim P+d_r-1}.
\]
In particular,
\[
[\Coh_{\mathbf d}(X)]
=\lim_{N\to\infty}[\Quot^X_{\widehat{\mathbf d}}(\mathcal O_X^N)]\L^{-Nd_r}.
\]
\end{lemma}
\begin{proof}
The vector bundle identity $[\mathcal V_{N,\mathbf d}]=\L^{Nd_r}[\Coh_{\mathbf d}(X)]$ and \eqref{eq:quot-stable-locus} give
\[
[\Coh_{\mathbf d}(X)]
-[\Quot^X_{\widehat{\mathbf d}}(\mathcal O_X^N)]\L^{-Nd_r}
=[\mathcal Z_{N,\mathbf d}]\L^{-Nd_r}.
\]
Apply Lemma~\ref{lem:unstable-locus}; the resulting filtration bound tends to $-\infty$ as $N\to\infty$.
\end{proof}

\begin{proof}[Proof of Theorem~\ref{thm:nestQC}]
Fix $\mathbf n\in\N^r$ and put $\mathbf d=\mathbf d(\mathbf n)$. The grading vector of $\widehat{\mathbf d}$ is $(n_r,\ldots,n_1)$, since
\[
\widehat d_i-\widehat d_{i-1}
=d_{r-i+1}-d_{r-i}=n_{r-i+1},\qquad \widehat d_0=0.
\]
As $d_r=n_1+\cdots+n_r$, Lemma~\ref{thm:A-effective} gives
\[
\begin{aligned}
[t_1^{n_1}\cdots t_r^{n_r}]
Z^X_{\mathcal O_X^N}(\L^{-N}t_r,\ldots,\L^{-N}t_1)
&=\L^{-Nd_r}[\Quot^X_{\widehat{\mathbf d}}(\mathcal O_X^N)]\\
&\longrightarrow[\Coh_{\mathbf d}(X)], \quad N\to \infty.
\end{aligned}
\]
This proves the claimed coefficientwise limit.
\end{proof}

\section{Flag zeta functions over a smooth curve}
Recall that Theorem~\ref{thm:flag-cl} follows directly from \cite{monavariricolfi} and Theorem~\ref{thm:nestQC}. For completeness, we give a short proof of the analogous result in the formal setting.
\begin{lemma}\label{lem:flag-cl}
    Let $r\geq 0$ and $R=k[\![t]\!]$. Then
 \begin{equation}
 {Z}^{R}_{R^d}(t_1,\dots,t_r)=\prod_{j=1}^r\prod_{i=0}^{d-1} \frac{1}{1-\L^{i} t_j},
\qquad
\Zhat_R(t_1,\dots,t_r)=\prod_{j=1}^r \prod_{i=1}^\infty \frac{1}{1-\L^{-i}t_j}.
    \end{equation}
\end{lemma}
\begin{proof}

We know that $ {Z}^{R}_{R^d}(t)=\prod_{i=0}^{d-1} \frac{1}{1-\L^{i} t}$. So it suffices to show that $ {Z}^{R}_{R^d}(t_1,\dots,t_r)=\prod_{i=1}^r {Z}^{R}_{R^d}(t_i)$. One can show it by choosing the steps inductively and noting that every torsion-free $R$-module of rank $d$ is isomorphic to $R^d$. To make this argument rigorous, we use the fact that the universal bundle over an Iwahoric cell (see \cite[\S 5]{huangjiang2023punctual}) is trivialized.

The second identity follows from the first by the Quot to Coh formula of Theorem~\ref{thm:nestQC}.
\end{proof}

Through this proof, we emphasize that the symmetry of the $\Zhat_R(t_1,\dots,t_r)$ and ${Z}^{R}_{R^d}(t_1,\dots,t_r)$ in the variables and their factorizability as a product of one-step zeta functions owe a lot to the fact that  the stack of rank $d$ torsion-free sheaves on $R$ has only one point, $R^d$, for the smooth germ $R=k[\![T]\!].$ Similar properties should not be expected when $R$ is singular.

\section{Flag zeta functions over a formal cusp}In this section, we prove Theorem~\ref{eq:nest cusp}(a)(b)(c)(d). Recall the notions from Definition~\ref{def:defcore}. Let $$ {Z}_{\mathcal{E}_a^{(d)};\mathcal{E}_b^{(d)}}^{R}(t):= \sum_{n\in \mathbb{N}}[\Quot_{n;\mathcal{E}_b^{(d)}}^R(\mathcal{E}_a^{(d)})]t^n.$$ This is an instance of the refined zeta function discussed in \S\ref{ref:2}. It is the $(a,b)$-th entry of $H^{(d)}_{R}(t)$. We also note that $$\Quot_{n}^R(\mathcal{E}_a^{(d)})=\bigsqcup_{b=0}^d \Quot_{n;\mathcal{E}_b^{(d)}}^R(\mathcal{E}_a^{(d)}),\;\;\;\;\; {Z}_{\mathcal{E}_a^{(d)}}^{R}(t)=\sum_{b=0}^d {Z}_{\mathcal{E}_a^{(d)};\mathcal{E}_b^{(d)}}^{R}(t).$$

\subsection{Computing $ {Z}_{M;\mathcal{E}_b}^{R}(t)$}\label{sub:computeZ^R} Fix $d$. In this section and the next, we will  omit the superscript ``$(d)$'' from $\mathcal{E}_a^{(d)}$ to ease notation. Let $M$ be a rank $d$ torsion-free $R$-module and let $\ut{M}\subseteq M$ be a rank $d$ $\tl{R}$-module. We have the following  locally closed subset of $\Quot_n^R(M)$, defined in \cite[\S 6]{huangjiang2023torsionfree}:
\begin{equation}
    B_R(M;\ut M;l;n) := \set*{
        L\in \Quot_n^R(M): \tl RL\supeq \ut M, \; \rk(\tl RL/\ut M)=l
    }.
\end{equation}
Here the $\rk(-)$ means the minimal number of generators as an $\tl{R}$-module. It follows from \cite[Theorem~6.2]{huangjiang2023torsionfree} that $$Z_M^R(t) = \sum_{l}  \frac{1}{(\L^lt;\L)_{d-l}} \sum_n [B_R(M;\ut M;l;n)] \, t^n.$$
\begin{remark}\label{rmk:structures1}
   The proof of this formula relies on the following structural result proved in \cite[\S 6]{huangjiang2023torsionfree}: 
   \begin{enumerate}
        \item For each $0\leq l\leq d$ and  $n\in \mathbb{N}$, there exists a collection of varieties  $\mathcal{V}_{l,n,\gamma},\gamma\in \Gamma_{l,n}$ where $\Gamma_{l,n}$ is a finite index set, such that we have a decomposition $$\Quot^R_n(M)=\bigsqcup_{l=0}^d\bigsqcup_{n_1+n_2=n}\bigsqcup_{\gamma\in \Gamma_{l,n_1}}\left(\mathcal{V}_{l,n_1,\gamma}\times B_R(M;\ut M;l;n_2)\right).$$
\item Let $x\in \mathcal{V}_{l,n_1,\gamma}\times B_R(M;\ut M;l;n_2)\subseteq \Quot^R_{n_1+n_2}(M)$ be a field valued point. Then the lattice determined by $x$ is isomorphic to the lattice determined by the projection of $x$ to $ B_R(M;\ut M;l;n_2)$.
      \item $\sum_{n\geq 0}\sum_{\gamma\in\Gamma_{l,n}}[\mathcal{V}_{l,n,\gamma}]t^{n}=\frac{1}{(\L^lt;\L)_{d-l}}$.
    \end{enumerate}
The technically demanding section \cite[\S 6]{huangjiang2023torsionfree} is essentially proving the existence of $\mathcal{V}_{l,n,\gamma}$ by  building
 them explicitly. But for our application, existence is enough.
\end{remark}

Now let $ B_{R;\mathcal{E}_b}(M;\ut M;l;n)\subseteq  B_R(M;\ut M;l;n)$ be the constructible subset that consists of points classifying modules that are isomorphic to $\mathcal{E}_b$. Then Remark~\ref{rmk:structures1} implies that
\begin{equation}\label{eq:refined Quot zeta}
   Z_{M;\mathcal{E}_b}^R(t) = \sum_{l}  \frac{1}{(\L^lt;\L)_{d-l}} \sum_n [B_{R;\mathcal{E}_b}(M;\ut M;l;n)] \, t^n.
\end{equation}
Now we set up our concrete computation when $M=\mathcal{E}_a$. To make everything explicit, we present $\mathcal{E}_a$ as a lattice in $\bigoplus_{i=1}^d k(\!(T)\!)e_i$ of the form
$\bigoplus_{i=1}^a T^2\widetilde{R} e_i\oplus \bigoplus_{j=a+1}^{d} Re_j $. Fix $\ut{\mathcal{E}_a}:= T^2\tl R^d=\bigoplus_{i=1}^d T^2\widetilde{R} e_i$.
\begin{lemma}\label{eq:B_R}
 Under the assumptions above, we have $$[B_{R;\mathcal{E}_b}(\mathcal{E}_a;\ut{\mathcal{E}_a};d-l;n)] =\left\{ \begin{aligned}
      &{\qbinom{d-a}{d-l}_{\L}\qbinom{l}{b}_{\L}}\L^{(l-b)d},\,\,&n=2l-a-b,\\
      &0,\,\,&\text{otherwise}.
  \end{aligned}\right.$$
In particular,
    $$ {Z}_{\mathcal{E}_a;\mathcal{E}_b}^{R}(t)=\sum_{l=\max\{a,b\}}^{d}{\qbinom{d-a}{d-l}_{\L}\qbinom{l}{b}_{\L}}\frac{\L^{(l-b)d}t^{2l-a-b}}{(\L^{d-l}t\,;\,\L)_{l}}.$$
\end{lemma}
We will use the observation that the isomorphism type of a rank-$d$ lattice $L$ is determined by $\dim_k(\tl RL/L)$. Indeed, the types $\mathcal E_b$ exhaust all possibilities, and
\[
\dim_k(\tl R\mathcal E_b/\mathcal E_b)=d-b,
\qquad L\simeq\mathcal E_b\ \Longleftrightarrow\ b=d-\dim_k(\tl RL/L).
\]
\begin{proof}
Set $U=\Span_k\{e_{a+1},\dots,e_d\}$, and let $\mathcal B$ be the reduced locus of lattices $L\subseteq\mathcal E_a$ with $\tl RL\supseteq T^2\tl R^d$. Write $\mathit{in}_i(L)$ for the space of vectors $v\in k^d$ occurring in elements $vT^i+O(T^{i+1})\in L$. For $L\subseteq\mathcal E_a$, the condition $\tl RL\supseteq T^2\tl R^d$ is equivalent to $\mathit{in}_2(L)=k^d$: the forward implication follows by comparing coefficients in $\tl RL=L+TL$ and using the absence of $T^1$-terms in $\mathcal E_a$; conversely, elements $e_iT^2+O(T^3)\in L$ form a $\tl R$-basis of $T^2\tl R^d$. Moreover, the conductor $T^2\tl R\subseteq R$ gives
\[
T^4\tl R^d\subseteq T^2\tl RL\subseteq L.
\]
Taking $\bbar L=L/T^4\tl R^d$, we may therefore regard $\mathcal B$ as the reduced locus parametrizing $k$-subspaces
\[
\bbar L\subseteq W:=\mathcal E_a/T^4\tl R^d=U\oplus k^dT^2\oplus k^dT^3,
\]
such that $R\bbar L\subseteq\bbar L$ and $\mathit{in}_2(\bbar L)=k^d$. Let $\mathcal B_{n_0,n_3}\subseteq \mathcal B$ be the reduced rank stratum where $\dim\mathit{in}_i(L)=n_i$ for $i=0,3$.

Let $\mathcal P_{n_0,n_3}$ parametrize quadruples $(V_0,V_3,\bbar\phi_0,\bbar\phi_2)$ with
\[
V_0\subseteq U,\quad V_0\subseteq V_3\subseteq k^d,\quad
\dim V_0=n_0,\quad\dim V_3=n_3,
\]
\[
\bbar\phi_0\in\Hom_k(V_0,k^d/V_3),\qquad
\bbar\phi_2\in\Hom_k(k^d,k^d/V_3),\quad\bbar\phi_2|_{V_0}=0.
\]
Its scheme structure is given by the tower
\[
\mathcal P_{n_0,n_3}\longrightarrow\mathcal Y_{n_0,n_3}
\longrightarrow\Gr(U,n_0),
\]
where $\mathcal Y_{n_0,n_3}=\Gr_{\Gr(U,n_0)}(\O^d/\mathcal S_0,n_3-n_0)$ and $\mathcal P_{n_0,n_3}$ is the total space of
\[
\mathcal Hom(\mathcal S_0,\O^d/\mathcal S_3)
\oplus\mathcal Hom(\O^d/\mathcal S_0,\O^d/\mathcal S_3)
\]
for the tautological flag $\mathcal S_0\subseteq\mathcal S_3$. In particular,
\[
[\mathcal P_{n_0,n_3}]
=\qbinom{d-a}{n_0}_{\L}\qbinom{d-n_0}{n_3-n_0}_{\L}\L^{d(d-n_3)}.
\]
We claim that $\mathcal P_{n_0,n_3}\simeq \mathcal B_{n_0,n_3}$.

Given a quadruple in $\mathcal P_{n_0,n_3}$, define
\[
\begin{aligned}
\bbar L={}&\{v-\phi_0(v)T^3:v\in V_0\}
\oplus\{vT^2-\phi_2(v)T^3:v\in k^d\}\oplus V_3T^3,
\end{aligned}
\]
using any lifts $\phi_0,\phi_2$ to $k^d$. The term $V_3T^3$ absorbs the choice of lifts. Equivalently, $\bbar L$ is the kernel of the surjective map
\[
\begin{aligned}
\theta:V_0\oplus k^dT^2\oplus k^dT^3&\longrightarrow k^d/V_3,\\
v_0+v_2T^2+v_3T^3&\longmapsto
\bbar\phi_0(v_0)+\bbar\phi_2(v_2)+(v_3\bmod V_3).
\end{aligned}
\]
The conditions $V_0\subseteq V_3$ and $\bbar\phi_2|_{V_0}=0$ make $\bbar L$ an $R$-submodule, and its degree-two generators give $\mathit{in}_2(\bbar L)=k^d$. Thus $\bbar L\in \mathcal B_{n_0,n_3}$. This defines a map $\Phi:\mathcal P_{n_0,n_3}\to \mathcal B_{n_0,n_3}$.

To see that it is a morphism into the ambient Grassmannian, pull the tautological bundles $\mathcal S_0,\mathcal S_3$ back to $\mathcal P_{n_0,n_3}$ and write $\O=\O_{\mathcal P_{n_0,n_3}}$. Using the universal maps $\bbar\phi_0,\bbar\phi_2$, define
\[
\begin{aligned}
\vartheta:\mathcal S_0\oplus\O^dT^2\oplus\O^dT^3&\longrightarrow\O^d/\mathcal S_3,\\
v_0+v_2T^2+v_3T^3&\longmapsto
\bbar\phi_0(v_0)+\bbar\phi_2(v_2)+(v_3\bmod\mathcal S_3).
\end{aligned}
\]
Since $\vartheta$ is surjective, its kernel commutes with base change. Thus the fiber of $\ker(\vartheta)\subseteq W\otimes_k\O$ over each quadruple is precisely $\bbar L\subseteq W$, of constant dimension $n_0+d+n_3$. Since $\mathcal P_{n_0,n_3}$ is reduced, this kernel is a subbundle, and the resulting Grassmannian morphism factors through $\mathcal B_{n_0,n_3}$. Hence $\Phi$ is a morphism.

Conversely, given $\bbar L\in\mathcal B_{n_0,n_3}$, recover $V_0=\mathit{in}_0(\bbar L)$ and $V_3=\mathit{in}_3(\bbar L)$. Equivalently,
\[ V_0=\im(\bbar L\to U),\qquad V_3T^3=\bbar L\cap k^d T^3,\]
where $\bbar L\to U$ is restricted from the projection map $W\to U$.

The inclusion of $k^dT^3$ induces the right-hand map in
\[
V_0\oplus k^dT^2
\longrightarrow
\frac{V_0\oplus k^dT^2\oplus k^dT^3}{\bbar L}
\xleftarrow[\ \sim\ ]{\alpha}
(k^d/V_3)T^3.
\]
This map is injective by the definition of $V_3$. It is surjective because, modulo $\bbar L$, we can eliminate the constant coefficient using the definition of $V_0$, and then the degree-two coefficient using $\mathit{in}_2(\bbar L)=k^d$. Composing the left-hand map with $\alpha^{-1}$ defines linear maps $\bbar\phi_0,\bbar\phi_2$ by
\[
v+uT^2\longmapsto
\bigl(\bbar\phi_0(v)+\bbar\phi_2(u)\bigr)T^3.
\]
For $v\in V_0$, multiplying an element of $\bbar L$ with constant coefficient $v$ by $T^3$ and $T^2$ gives $vT^3,vT^2\in\bbar L$. Hence $V_0\subseteq V_3$ and $\bbar\phi_2|_{V_0}=0$. The recovered quadruple belongs to $\mathcal P_{n_0,n_3}$, and its lattice is $\bbar L$, since the remainder of $v+uT^2+wT^3$ is
$\bigl(\bbar\phi_0(v)+\bbar\phi_2(u)+(w\bmod V_3)\bigr)T^3$. This completes the set-theoretical construction of the inverse map of $\Phi$. 

To see $\Phi^{-1}$ is a morphism, we use a similar argument again, noting that all linear maps used in the construction can be upgraded into maps between vector bundles, and the reducedness of $\mathcal{B}_{n_0,n_3}$ and the constant fiber dimension argument account for every required subbundle claim. This establishes the required isomorphism $\Phi$.

It remains to identify the required stratum. For the lattice above,
\[
\tl RL=V_0\oplus TV_0\oplus T^2\tl R^d,
\]
and comparison modulo $T^4$ gives
\begin{align*}
\rk(\tl RL/\ut{\mathcal E_a})&=n_0,\\
\dim_k(\mathcal E_a/L)&=2d-a-n_0-n_3,\\
\dim_k(\tl RL/L)&=d+n_0-n_3.
\end{align*}
By the observation preceding the proof, the required conditions are therefore
\[
n_0=d-l,\qquad n_3=d+b-l,\qquad n=2l-a-b.
\]
Therefore the required locus is $\mathcal B_{d-l,d+b-l}$ when $n=2l-a-b$, and is empty otherwise; this is an equality of reduced subvarieties in the same Quot scheme. Substituting in the class of $\mathcal P_{n_0,n_3}$ proves the first formula; \eqref{eq:refined Quot zeta} gives the second.
\end{proof}

\subsection{Computing ${Z}^{R}_{\mathcal{E}_a}(t_1,\dots,t_r)$}
 \begin{lemma}\label{eq:reductionformula}
$ {Z}^{R}_{\mathcal{E}_a}(t_1,\dots,t_r)=\sum_{b=0}^d {Z}^{R}_{\mathcal{E}_a ; \mathcal{E}_b}(t_1)\cdot {Z}^{R}_{\mathcal{E}_b}(t_2,\dots,t_{r})$. Iterating, we get \begin{align*}
    {Z}^{R}_{\mathcal{E}_a}(t_1,\dots,t_r)&=\sum_{0\leq a_1,\dots,a_r\leq d} {Z}^{R}_{\mathcal{E}_{a} ;\mathcal{E}_{a_
1}}(t_1) {Z}^{R}_{\mathcal{E}_{a_1} ;  \mathcal{E}_{a_2}}(t_{2})\cdots  {Z}^R_{\mathcal{E}_{a_{r-1}}; \mathcal{E}_{a_{r}}}(t_{r}).
\end{align*}
 \end{lemma}
\begin{proof}
First consider $S=\Quot^R_{n;\mathcal E_b}(\mathcal E_a)$, with universal lattice $\mathcal L$. We show that $S$ admits a finite stratification on which $\mathcal L$ is isomorphic to the pullback of $\mathcal E_b$. Write $R_S=\O_S[\![T^2,T^3]\!]$ and $\tl R_S=\O_S[\![T]\!]$. Since $S$ is reduced and the universal quotient is locally free of rank $n$ over $S$, the fiberwise identity $\m^n(\mathcal E_{a,S}/\mathcal L)=0$ holds over $S$. Hence $T^{2n+2}\tl R_S^d\subseteq\mathcal L$. Set $\tl{\mathcal L}=\mathcal L+T\mathcal L$; the conductor gives $T^2\tl{\mathcal L}\subseteq\mathcal L$.

Work in the finite quotient $\tl R_S^d/T^{2n+4}\tl R_S^d$. Take successive rank strata for the maps defining $\tl{\mathcal L}$, $T\tl{\mathcal L}$, $T^2\tl{\mathcal L}$, and $\mathcal L\cap T\tl{\mathcal L}$, so that their subquotients are vector bundles and commute with base change. On each stratum there is an exact sequence of vector bundles
\[
0\longrightarrow\mathcal W:=\frac{\mathcal L\cap T\tl{\mathcal L}}{T^2\tl{\mathcal L}}
\longrightarrow\frac{\mathcal L}{T^2\tl{\mathcal L}}
\longrightarrow\frac{\tl{\mathcal L}}{T\tl{\mathcal L}}\longrightarrow0
\]
of ranks $b,d+b,d$. Multiplication by $T$ identifies $\tl{\mathcal L}/T\tl{\mathcal L}$ with $T\tl{\mathcal L}/T^2\tl{\mathcal L}$, so we can choose a local basis of $\tl{\mathcal L}/T\tl{\mathcal L}$ adapted to $\mathcal W$ and lift it to $e_1,\dots,e_d\in\mathcal L$. Successive lifting in $T$ shows that these lifts form a $\tl R_S$-basis of $\tl{\mathcal L}$, and
\[
\mathcal L=\bigoplus_{i=1}^b\tl R_S e_i\ \oplus\!\!\bigoplus_{i=b+1}^d R_S e_i
\simeq\mathcal E_{b,S},
\]
where $\mathcal E_{b,S}$ denotes the pullback of $\mathcal E_b$. Refining finite open covers to a locally closed stratification proves the claim.

Now fix $\mathbf d=(d_1,\dots,d_r)$ and $\mathbf d'=(d_2-d_1,\dots,d_r-d_1)$. The map
$\Quot^R_{\mathbf d}(\mathcal E_a)\to\Quot^R_{d_1}(\mathcal E_a)$ sends $(\mathcal{K}_1\supeq \dots \supeq \mathcal{K}_r)$ to $\mathcal{K}_1$. Over each stratum $S_\alpha\subseteq\Quot^R_{d_1;\mathcal E_b}(\mathcal E_a)$ constructed above, the relative flag Quot functor identifies its inverse image with $S_\alpha\times\Quot^R_{\mathbf d'}(\mathcal E_b)$. Additivity in the Grothendieck ring therefore gives
\[
[\Quot^R_{\mathbf d}(\mathcal E_a)]
=\sum_{b=0}^d[\Quot^R_{d_1;\mathcal E_b}(\mathcal E_a)]
[\Quot^R_{\mathbf d'}(\mathcal E_b)].
\]
Summing over $\mathbf{d}$ gives the required formula.
\end{proof}

\begin{proof}[Proof of Theorem~\ref{eq:nest cusp}(a)(b)(c)(d)]
    Combining Lemma~\ref{eq:B_R} and Lemma~\ref{eq:reductionformula}, parts (a)(c) of Theorem~\ref{eq:nest cusp} follow. Then apply Theorem~\ref{thm:nestQC} to get (b)(d).
\end{proof}

\section{One-step normalized zeta functions}
\label{sec:one-step-interpolation}

In this section we extract and simplify the one-step consequence of
Theorem~\ref{eq:nest cusp}.  The simplification rests on the following finite
Gaussian-binomial identity. Here $q$ is an indeterminate. We use the Rogers--Szeg\H{o} polynomial
\[
 h_m(x;q):=\sum_{u=0}^{m}\qbinom{m}{u}_{q}x^u.
\]

\begin{lemma}\label{lem:finite-gaussian-collapse}
Let $a,n\geq 0$ and put $d=a+n$.  Then
\begin{align}
 &\sum_{r=0}^{n}
 \qbinom{n}{r}_{q}(z;q)_r z^{n-r}
 h_{d-r}(q^d z;q) \notag\\
 &\qquad=
 \sum_{i=0}^{a}\sum_{j=0}^{n}
 q^{d(i+j)+(a-i)j}
 \qbinom{a}{i}_{q}\qbinom{n}{j}_{q}
 z^{i+2j}.
 \label{eq:finite-gaussian-collapse}
\end{align}
\end{lemma}

\begin{proof}
By the finite $q$-binomial theorem,
\[
 (z;q)_r
 =\sum_{h=0}^{r}(-1)^h q^{\binom h2}
   \qbinom{r}{h}_{q}z^h.
\]
After also expanding $h_{d-r}(q^dz;q)$, the left-hand side of
\eqref{eq:finite-gaussian-collapse} becomes
\[
 \sum_{r=0}^{n}\sum_{h=0}^{r}\sum_{p=0}^{d-r}
 (-1)^h q^{\binom h2+dp}
 \qbinom{n}{r}_{q}
 \qbinom{r}{h}_{q}
 \qbinom{d-r}{p}_{q}
 z^{n-r+h+p}.
\]
Set
\[
 s=r-h,\qquad j=n-s,\qquad i=p-j.
\]
Then $0\leq j\leq n$, while $i$ initially ranges over a finite subset of
$\mathbb Z$.  Using
\[
 \qbinom{n}{s+h}_{q}\qbinom{s+h}{h}_{q}
 =\qbinom{n}{j}_{q}\qbinom{j}{h}_{q},
\]
the preceding triple sum can be reorganized as
\begin{align*}
 \sum_{j=0}^{n}\sum_{i\in\mathbb Z}
 q^{d(i+j)}z^{i+2j}\qbinom{n}{j}_{q}
 \sum_{h=0}^{j}(-1)^h q^{\binom h2}
 \qbinom{j}{h}_{q}
 \qbinom{a+j-h}{i+j}_{q}.
\end{align*}
It remains to use
\begin{equation}\label{eq:mobius-gaussian}
 \sum_{h=0}^{j}(-1)^h q^{\binom h2}
 \qbinom{j}{h}_{q}
 \qbinom{a+j-h}{i+j}_{q}
 =q^{(a-i)j}\qbinom{a}{i}_{q}.
\end{equation}

Here and below Gaussian binomial coefficients with an inadmissible lower
index are zero.  For completeness, identity \eqref{eq:mobius-gaussian} has a direct
finite-geometric proof. Temporarily let $q$ be a prime power. Let $A$ and $J$ be vector spaces of dimensions $a$
and $j$ over $\mathbb F_q$.  We count $(i+j)$-dimensional subspaces
$W\subseteq A\oplus J$ whose projection to $J$ is surjective.  Directly, one
first chooses the $i$-dimensional subspace $W\cap A$ and then a linear map
$J\to A/(W\cap A)$, giving
\[
 q^{(a-i)j}\qbinom{a}{i}_{q}
\]
possibilities.  Alternatively, apply M\"obius inversion to the lattice of
subspaces of $J$.  A codimension-$h$ subspace of $J$ has M\"obius coefficient
$(-1)^h q^{\binom h2}$, and there are $\qbinom{j}{h}_{q}$ such subspaces.  The
number of $(i+j)$-dimensional subspaces contained in the corresponding
$(a+j-h)$-dimensional space is $\qbinom{a+j-h}{i+j}_{q}$.  This gives the
left-hand side of \eqref{eq:mobius-gaussian}.  This argument also gives zero
when $i<0$ or $i>a$, since no such $W$ exists.  Since the resulting identity
holds for every prime power $q$, it is an identity in $\mathbb Z[q]$.
Substitution into the preceding triple sum proves the lemma.
\end{proof}

We now apply the lemma to the cusp.  Recall that
\[
 \mathcal E_a^{(d)}=\widetilde R^{a}\oplus R^{d-a},
 \qquad
 N^R_{\mathcal E_a^{(d)}}(t)
 =(t;\L)_d Z^R_{\mathcal E_a^{(d)}}(t).
\]

\begin{theorem}\label{thm:one-step-interpolation}
For every $0\leq a\leq d$, one has
\begin{equation}\label{eq:one-step-interpolation}
 N^R_{\mathcal E_a^{(d)}}(t)
 =
 \sum_{i=0}^{a}\sum_{j=0}^{d-a}
 \L^{d(i+j)+(a-i)j}
 \qbinom{a}{i}_{\L}
 \qbinom{d-a}{j}_{\L}
 t^{i+2j}.
\end{equation}
Equivalently,
\begin{equation}\label{eq:one-step-interpolation-m}
 N^R_{\mathcal E_a^{(d)}}(t)
 =\sum_{m=0}^{d}\L^{dm}
 \sum_i
 \L^{(a-i)(m-i)}
 \qbinom{a}{i}_{\L}
 \qbinom{d-a}{m-i}_{\L}
 t^{2m-i},
\end{equation}
where a Gaussian binomial coefficient is understood to vanish when its lower
index is outside its natural range.
\end{theorem}

\begin{proof}
By Theorem~\ref{eq:nest cusp}(a) and (c),
\begin{align*}
 Z^R_{\mathcal E_a^{(d)}}(t)
 &=\sum_{b=0}^{d}\sum_{l=\max\{a,b\}}^{d}
 \qbinom{d-a}{d-l}_{\L}\qbinom{l}{b}_{\L}
 \frac{\L^{(l-b)d}t^{2l-a-b}}
 {(\L^{d-l}t;\L)_l}.
\end{align*}
For fixed $l$, set $u=l-b$.  The sum over $b$ is
\[
 t^{l-a}\sum_{u=0}^{l}
 \qbinom{l}{u}_{\L}(\L^dt)^u
 =t^{l-a}h_l(\L^dt;\L).
\]
Moreover,
\[
 \frac{(t;\L)_d}{(\L^{d-l}t;\L)_l}=(t;\L)_{d-l}.
\]
Consequently, setting $r=d-l$ gives
\begin{equation}\label{eq:one-step-before-collapse}
 N^R_{\mathcal E_a^{(d)}}(t)
 =\sum_{r=0}^{d-a}
 \qbinom{d-a}{r}_{\L}
 (t;\L)_r t^{d-a-r}
 h_{d-r}(\L^dt;\L).
\end{equation}
Apply Lemma~\ref{lem:finite-gaussian-collapse} with
$q=\L$, $z=t$, and $n=d-a$.  This gives
\eqref{eq:one-step-interpolation}.  Formula
\eqref{eq:one-step-interpolation-m} follows by setting $m=i+j$.
\end{proof}

There is a useful geometric interpretation of the positive summands in
Theorem~\ref{thm:one-step-interpolation}.  Fix an $a$-dimensional subspace
$A\subseteq k^d$.  For $0\leq i\leq m\leq d$, let
\[
 \Gr(m,d)_{A,i}
 :=\{W\in\Gr(m,d):\dim(W\cap A)=i\}.
\]
The map
\[
 \Gr(m,d)_{A,i}\longrightarrow
 \Gr(i,A)\times\Gr(m-i,k^d/A),
 \qquad
 W\longmapsto (W\cap A,\operatorname{im}(W\to k^d/A)),
\]
is a vector bundle of rank $(a-i)(m-i)$.  Hence
\begin{equation}\label{eq:grassmannian-stratum}
 [\Gr(m,d)_{A,i}]
 =\L^{(a-i)(m-i)}
 \qbinom{a}{i}_{\L}
 \qbinom{d-a}{m-i}_{\L}.
\end{equation}
Formula \eqref{eq:one-step-interpolation-m} may therefore be written as the
motivic weighted sum
\begin{equation}\label{eq:one-step-grassmannian}
 N^R_{\mathcal E_a^{(d)}}(t)
 =\sum_{m=0}^{d}\L^{dm}
 \sum_{i}[\Gr(m,d)_{A,i}]t^{2m-i}.
\end{equation}
Thus the exponent of $t$ records both the dimension of $W$ and its relative
position with respect to the fixed $a$-plane $A$.

\begin{corollary}\label{cor:one-step-endpoints}
The two extreme cases of Theorem~\ref{thm:one-step-interpolation} are
\begin{align}
 N^R_{R^d}(t)
 &=\sum_{j=0}^{d}
 \L^{dj}\qbinom{d}{j}_{\L}t^{2j},
 \label{eq:endpoint-R}\\
 N^R_{\widetilde R^d}(t)
 &=\sum_{i=0}^{d}
 \L^{di}\qbinom{d}{i}_{\L}t^i.
 \label{eq:endpoint-Rtilde}
\end{align}
In particular,
\[
 N^R_{R^d}(t)=N^R_{\widetilde R^d}(t^2).
\]
These recover the $R_{2,3}$ cases of
\cite[Proposition~8.5 and Theorem~9.1]{huangjiang2023torsionfree}.
\end{corollary}

\begin{proof}
Set $a=0$ and $a=d$, respectively, in
\eqref{eq:one-step-interpolation}.
\end{proof}

We finish by placing the interpolation in its Rogers--Ramanujan form.  Put
$\mathbb q=\L^{-1}$.  The relation
\[
 \qbinom{n}{r}_{\L}
 =\L^{r(n-r)}\qbinom{n}{r}_{\mathbb q}
\]
transforms \eqref{eq:one-step-interpolation} into
\begin{equation}\label{eq:one-step-q-series}
 N^R_{\mathcal E_a^{(d)}}(\L^{-d}t)
 =\sum_{i=0}^{a}\sum_{j=0}^{d-a}
 \mathbb q^{i^2-ai+ij+j^2}
 \qbinom{a}{i}_{\mathbb q}
 \qbinom{d-a}{j}_{\mathbb q}
 t^{i+2j},
\end{equation}
which converges as $d\to \infty$; cf.~\eqref{eq:intro-a2-slice}.

As two elementary observations, setting $t=1$ in
\eqref{eq:one-step-grassmannian} and forgetting the relative-position
stratification gives
\begin{equation}\label{eq:one-step-t-one}
 N^R_{\mathcal E_a^{(d)}}(1)
 =\sum_{m=0}^{d}\L^{dm}\qbinom{d}{m}_{\L},
\end{equation}
which is independent of $a$, consistent with a general geometric fact \cite[Corollary~6.4]{huangjiang2023torsionfree}.  Moreover,
\[
 \deg_t N^R_{\mathcal E_a^{(d)}}(t)=2d-a,
 \qquad
 [t^{2d-a}]N^R_{\mathcal E_a^{(d)}}(t)=\L^{d^2}.
\]

\begin{example}\label{ex:rank-two-interpolation}
For $d=2$, Theorem~\ref{thm:one-step-interpolation} gives
\begin{align*}
 N^R_{R^2}(t)
 &=1+(\L^2+\L^3)t^2+\L^4t^4,\\
 N^R_{\widetilde R\oplus R}(t)
 &=1+\L^2t+\L^3t^2+\L^4t^3,\\
 N^R_{\widetilde R^2}(t)
 &=1+(\L^2+\L^3)t+\L^4t^2.
\end{align*}
\end{example}

\section{The stable Hecke algebra}
In this section we prove Theorem~\ref{thm:commuteHecke} and deduce Theorem~\ref{eq:nest cusp}(e)(f) as a consequence.

We set by default that $\widehat{H}_{R,-1}=0$. Note also that $\widehat{H}_{R,0}=I$, the infinite rank matrix with 1's on the diagonal. An easy computation using (\ref{eq:limit QUOT}) shows that  $\widehat{H}_{R,1}$ is the matrix with $(a,b)$-th entry
\begin{equation}\label{eq:stable-hecke-one}
(\widehat{H}_{R,1})_{a,b}=\left\{\begin{aligned}
&\mathbb{q}[a]_\mathbb{q}   &b=a,  \\
&\mathbb{q}^{a+1}/(1-\mathbb{q}) &b=a+1,\\
&\mathbb{q}^{1-a}[a]_\mathbb{q}  &b=a-1,\\
&0 & \text{otherwise}.
\end{aligned} \right.
\end{equation}

By \eqref{eq:limit QUOT}, $(\widehat H_{R,m})_{a,b}=0$ for $|a-b|>m$. Thus products of these infinite matrices and their generating series, as well as the row sums in Theorem~\ref{eq:nest cusp}(d), involve only finitely many entries in each fixed degree. In particular, polynomials in the tridiagonal matrix $\widehat H_{R,1}$ are well-defined. Infinite products are interpreted entrywise, $\mathbb q$-adically in each coefficient of $t$; since $\mathbb q=\L^{-1}$, this also gives convergence in the dimension filtration (Definition~\ref{def:k0stck}).

\begin{theorem}\label{thm:recursion}
    The followings are true:\begin{enumerate}
         \item For $m\geq 0$, we have
\begin{equation}\label{eq:coreeq}
\widehat{H}_{R,1}\cdot\widehat{H}_{R,m}
 =[m+1]_\L\widehat{H}_{R,m+1}-[m]_\L\widehat{H}_{R,m}
 -\frac \L{\L-1}\widehat{H}_{R,m-1}.
\end{equation}
   \item Let $\mathbb{q}=\L^{-1}$. Then we have $\widehat{H}_R(0)=I$ and
\begin{equation}\label{eq:recure2}
    \widehat{H}_{R}(t)
 =\frac{1-\mathbb{q}t+\mathbb{q}t^2+(1-\mathbb{q})t\, \widehat{H}_{R,1}}{1-\mathbb{q}t}\, \widehat{H}_{R}(\mathbb{q}t).
\end{equation}
    \end{enumerate}
\end{theorem}
\begin{proof}
We prove (b); direct expansion implies (a). Write $q=\mathbb q$. For $l\in \Z$, consider the matrix $H_l(t)$ with $(a,b)$-entry
\[
H_{l;a,b}(t)=\qbinom{l}{b}_q
\frac{q^{l(l-a)-b(l-b)}t^{2l-a-b}}
{(q;q)_{l-a}(qt;q)_l}
\qquad(0\leq a,b\leq l), 
\]
and set $H_{l;a,b}(t)=0$ if $l<\max\{a,b\}.$

Then the formula \eqref{eq:limit QUOT} simply reads $\widehat H_R(t)=\sum_{l\in \Z} H_l(t)$. Put $\mathbf h_{l;b}(t)=H_{l}(t)e_b$, the $b$-th column of $H_l(t)$. 

Write $M(t)=(1-qt+qt^2)I+(1-q)t\widehat H_{R,1}$, then \eqref{eq:recure2} is equivalent to 
\begin{equation}\label{eq:recure3}
M(t)\widehat H_R(qt)e_b=(1-qt)\widehat H_R(t)e_b, \qquad b\ge 0.
\end{equation}

The key observation is
\begin{equation}\label{eq:hecke-summand-shift}
M(t)\mathbf h_{l;b}(qt)=(1-qt)\bigl(q^{l-b}\mathbf h_{l;b}(t)+(1-q^{l+1-b})\mathbf h_{l+1;b}(t)\bigr), \qquad b\ge 0, l\in \Z
\end{equation}
Summing over $l\in \Z$, telescoping gives \eqref{eq:recure3}, finishing the proof.

All sums above are coefficientwise finite in $t$. To verify
\eqref{eq:hecke-summand-shift}, its $a$-th component, expanded using
\eqref{eq:stable-hecke-one}, is
\begin{equation}\label{eq:entrywise}
    \begin{aligned}
&(1-q^{a+1}t+qt^2)H_{l;a,b}(qt)
 +q^{a+1}tH_{l;a+1,b}(qt)\\
&\qquad
 +q^{1-a}(1-q^a)tH_{l;a-1,b}(qt)\\
&\quad=(1-qt)\bigl(
 q^{l-b}H_{l;a,b}(t)
 +(1-q^{l+1-b})H_{l+1;a,b}(t)
 \bigr),
\end{aligned}
\end{equation}
where the third term on the left is omitted when $a=0$.

If $l<b-1$ or $l<a-1$, both sides vanish. If $l=b-1$, the only surviving term is the second term on the right, which has
coefficient $1-q^0=0$. If $l=a-1$, only the third term on the left and the second term on the right survive. If $l=a$, the second term on the left vanishes. We may therefore assume $l\geq b$ and
$l\ge a-1$, and put $r=l-a+1\geq0$.  
Removing the common factor
\[
\qbinom{l}{b}_q
\frac{q^{lr-b(l+1-b)-1}t^{2l-a-b}}
{(q;q)_r(q^2t;q)_l}
\]
reduces the desired equality \eqref{eq:entrywise}, term by term, to
\begin{align*}
&q^r(1-q^r)(1-q^{a+1}t+qt^2)
 +(q-q^r)(1-q^r)
 +q^{2r+1}(1-q^a)t^2\\
&\qquad=
 q(1-q^r)(1-q^{l+1}t)
 +q^{r+1}(1-q^{l+1})t^2;
\end{align*}
note that omitted terms at $a=0$ and $r=0,1$ still vanish in this polynomial identity. Direct expansion using $a+r=l+1$ verifies the polynomial identity, completing the proof.
\end{proof}
\begin{proof}[Proof of Theorem~\ref{thm:commuteHecke} and Theorem~\ref{eq:nest cusp}(e)(f)]
Define polynomials $P_m(x)\in \mathbb{Q}(\L)[x]$ recursively by $$[m+1]_{\L}P_{m+1}(x)
 =(x+[m]_{\L})P_m(x)+\frac \L{\L-1}P_{m-1}(x),\quad P_0(x)=1,P_1(x)=x.$$
 Then (\ref{eq:coreeq}) implies that $\widehat{H}_{R,m}=P_{m}(\widehat{H}_{R,1})$. Being polynomials in one matrix, the $\widehat H_{R,m}$ commute pairwise, proving Theorem~\ref{thm:commuteHecke}. Hence the series $\widehat H_R(t_i)$ commute, and permuting their factors in Theorem~\ref{eq:nest cusp}(d) proves the symmetry in part (f).
 
Finally, the equation \eqref{eq:recure2} has a unique solution with constant term $I$, given by the equivalent recursion \eqref{eq:coreeq}. Consider the product in Theorem~\ref{eq:nest cusp}(e), having constant term $I$. Its $s$-th factor is $I+O(\mathbb q^s)I+O(\mathbb q^s)\widehat H_{R,1}$. Though $\widehat H_{R,1}$ has unbounded negative powers of $\mathbb q$, we observe from the expression of the $s$-th factor that in any $t^m$-coefficient of the product, only powers $\widehat H_{R,1}^{\,j}$ with $0\le j\le m$ occur, contributing to bounded negative powers of $\mathbb q$ in any given entry. Thus, only finitely many terms in the infinite product contribute to a prescribed power of $\mathbb q$ in the $t^m$-coefficient of a given entry, so the product converges. The scalar prefactor $(\mathbb qt;\mathbb q)_\infty^{-1}$ converges by the same coefficientwise criterion. The full expression in Theorem~\ref{eq:nest cusp}(e) satisfies \eqref{eq:recure2} by separating the first factor from both sides. By uniqueness, it equals $\widehat H_R(t)$, proving Theorem~\ref{eq:nest cusp}(e).
\end{proof}

\appendix
\section{Explicit stratification of full flag $\Coh$ over a formal cusp} Let $R=k[\![x,y]\!]/(x^2-y^3)$, and let $\m$ be its maximal ideal. Let $\mathbf{d}=(1,\dots,d)$. We determine the motive of the flag coherent stack $\mathrm{Coh}_{\mathbf{d}}(R)$. Consider \begin{equation}\label{eq:upper_triangular}
V_d:=\set*{(X,Y)\in \Mat_d(k)^2\ \vrule \ \begin{array}{l}\text{$X,Y$ are strictly upper triangular}\\X^2=Y^3, XY=YX\end{array}}.\end{equation}
Note that $V_d$ and $\mathrm{C}_{\mathbf{d}}(R)$ have the same class in $\KVar{k}$  by Remark~\ref{rmk:indformal}. Now Corollary~\ref{cor:artin} implies that $$\frac{[V_d]}{(\L-1)^{d}\L^{d(d-1)/2}}=[\Coh_\mathbf{d}(R)]=[t_1t_2\cdots t_d]\Zhat_R(t_1,\dots,t_d),$$
so the knowledge of $[V_d]$ is a slice of Theorem~\ref{eq:nest cusp}(d)(e). We record below an independent computation of $[V_d]$ in \cite{huangjiang2023punctual}. We leave it to the reader to verify that the motivic implication of Theorem~\ref{cor:staircase_formula} is consistent with Theorem~\ref{eq:nest cusp}. A likely strategy is outlined in Remark~\ref{rmk:borel}.

\begin{theorem}\label{cor:staircase_formula}
There is a stratification
$$ V_d=\bigsqcup_{\substack{a+b=2d\\a\geq b\geq 0}} V\sbp{a,b}$$
for certain locally closed subsets $V\sbp{a,b}$. The motive of $V\sbp{a,b}$ can be computed inductively by
\begin{equation}
\begin{gathered}
[V\sbp{0,0}]=1, [V\sbp{a,b}]=0\text{ unless }a\geq b\geq 0, a\equiv b\pmod 2;\\
[V\sbp{a,b}]=\L^b [V\sbp{a-2,b}]+ (\L^{\frac{a+b-2}{2}}-\L^{b-1}) [V\sbp{a-1,b-1}] + (\L^a - \L^{\frac{a+b-2}{2}}) [V\sbp{a,b-2}].
\end{gathered}
\end{equation}

\end{theorem}

\subsection{Induction setup for the motive of $V_d$}
We use induction on $d$. Let $(X_{d+1},Y_{d+1})\in \Mat_{d+1}(k)^2$ be given by
\begin{equation}\label{eq:induction_rule}
X_{d+1} = \begin{bmatrix}
X_d & z \\
0 & 0
\end{bmatrix}, \quad
Y_{d+1} = \begin{bmatrix}
Y_d & w \\
0 & 0
\end{bmatrix},
\end{equation}
where $X_d, Y_d\in \Mat_d(k)$ and $z,w\in \Mat_{d\times 1}(k)$. By inspecting \eqref{eq:upper_triangular}, the pair $(X_{d+1},Y_{d+1})$ is in $V_{d+1}$ if and only if
\begin{equation}\label{eq:induction}
\begin{aligned}
(X_d,Y_d)&\in V_d; \\
X_d \,z &= Y_d ^2 \,w; \\
X_d \,w &= Y_d \,z.
\end{aligned}
\end{equation}

Thus, we have a morphism $\Phi_d: V_{d+1}\to V_d$ by sending $(X_{d+1},Y_{d+1})$ to $(X_d,Y_d)$, whose fibers are solution spaces of linear equations. Define the constructible function $a: V_d\to \N$ given by $a(X_d,Y_d):=\dim \Phi_d^{-1}(X_d,Y_d)$. In order to compute the motive of $V_{d+1}$, we need to compute the motive of the stratum $a^{-1}(i)$ for each $i\in \N$. Furthermore, in order to carry out the induction, we also need to understand how $a(X_{d+1},Y_{d+1})$ depends on the choice of $(X_d,Y_d)$ and $(z,w)$. This is a challenging task in general, but it turns out that we are able to reinterpret \eqref{eq:induction} as equations on modules over a ring, so that we can exploit tools in homological algebra.

\subsection{Some homological algebra}
Consider the $R$-module $M_d$ whose underlying vector space is $k^d$ and the multiplications by $x$ and $y$ are given by the matrices $X_d$ and $Y_d$, respectively. The condition $(X_d,Y_d)\in V_d$ ensures that the module structure of $M_d$ is well-defined. The module $M_d$ will be where $z$ and $w$ live. Consider the matrix
\begin{equation}
A=\begin{bmatrix}
x & -y^2\\
-y & x
\end{bmatrix}\in \Mat_2(R).
\end{equation}

Then \eqref{eq:induction} can be restated as $(z,w)\in \ker A_{M_d}$, where $A_{M_d}\in \End_R(M_d\oplus M_d)$ is the tensor product of $A:R^2\to R^2$ with $M_d$. We write $a(M_d):=a(X_d,Y_d)=\dim_k \ker A_{M_d}$. We also write $b(M_d)=\dim_k \im A_{M_d}$. By the rank-nullity theorem for $A_{M_d}$, we have $a(M_d)+b(M_d)=2d$ and $a(M_d)=\dim_k \coker A_{M_d}$.

Recall that we need to understand how $a(M_d)$ varies as $M_d=(X_d,Y_d)\in V_d$ varies. A general principle suggests that we shall fix a minimal resolution of $\coker A$:
\begin{equation}\label{eq:resolution}
\cdots \map[A] R^2 \map[A'] R^2  \map[A] R^2 \map \coker A \map 0,
\end{equation}
where
\begin{equation}
A'=\begin{bmatrix}
x & y^2\\
y & x
\end{bmatrix}.
\end{equation}

One can immediately verify $AA'=A'A=0$, so the above is indeed a chain complex.

\begin{remark}
Incidentally, $\coker A$ is isomorphic to $\m$, so $\coker A_{M_d} \cong \m \otimes_R M_d$, so that $a(M_d)=\dim_k \parens*{\m\otimes_R M_d}$. However, we are not able to make use of this observation. In fact, the rest of the proof below does not even use the fact that \eqref{eq:resolution} is exact, but only that it is a chain complex.
\end{remark}

Consider the $2$-periodic chain complex
\begin{equation}
C^\bullet: \cdots \map[A] R^2 \map[A'] R^2 \map[A] R^2 \map[A'] \cdots,
\end{equation}
where the degree convention is made such that a part of $C^\bullet$ reads $C^{-1} \map[A'] C^0 \map[A] C^1$. We explore some properties that hold \emph{universally} for $C^\bullet$, i.e., hold for the chain complex $C_M^\bullet:=C^\bullet \otimes_R M$ for any $R$-module $M$. Define the matrices
\begin{equation}
T=\begin{bmatrix}
0 & y \\
1 & 0
\end{bmatrix}, \quad
H=\begin{bmatrix}
0 & 1 \\
0 & 0
\end{bmatrix}, \quad H'=-H
\end{equation}
in $\Mat_2(R)$ and the matrix
\begin{equation}
K=\begin{bmatrix}
0 & 1 & 0 & -y \\
0 & 0 & 1 & 0
\end{bmatrix}
\end{equation}
in $\Mat_{2\times 4}(R)$.

The following matrix identities on $\Mat_2(R)$ or $\Mat_{2\times 4}(R)$ can be verified directly using the relation $x^2=y^3$ in $R$.
\begin{lemma}\label{lem:auxiliary_identities}
~
\begin{enumerate}
\item $[A,T]=[A',T]=0$.
\item $H'A+A'H=HA'+AH'=T^2$.
\item $
AK=H\begin{bmatrix}
A & O_{2\times 2}
\end{bmatrix}-T\begin{bmatrix}
T & -A'
\end{bmatrix}
$.
\item $
\begin{bmatrix}
 I_{2\times 2} & O_{2\times 2}
\end{bmatrix}
-TK = H\begin{bmatrix}
T & -A'
\end{bmatrix} + A'\begin{bmatrix}
O_{2\times 2} & H
\end{bmatrix}
$.
\item Let $g=\mathrm{diag}(1,-1)\in \Mat_2(R)$, and note that $g=g^{-1}$. Then $gAg=A'$ and $gTg=-T$.
\end{enumerate}
\end{lemma}

From (a), the matrix $T$ induces a chain map $T:C^\bullet\to C^\bullet$. In particular, for any $R$-module $M$, we have the induced chain map $T_M:C_M^\bullet\to C_M^\bullet$ and the induced maps on the cohomology of $C_M^\bullet$:
\begin{equation}
T^0_M: H^0(C_M^\bullet)\to H^0(C_M^\bullet), \quad T^1_M: H^1(C_M^\bullet)\to H^1(C_M^\bullet).
\end{equation}

Lemma \ref{lem:auxiliary_identities}(e) implies that the chain $C_M^\bullet$ together with the chain map $T_M$ ``looks the same'' in every degree:

\begin{lemma}\label{lem:iso_to_shift}
The matrix $g:C^i=R^2\to C^{i+1}=R^2$ induces a chain isomorphism $g: C^\bullet \to C^\bullet[1]$, where we use the unsigned shift: $(C^\bullet[1])^i=C^{i+1}$ and $d_{C[1]}^i=d_C^{i+1}$. Moreover, we have a commutative diagram of chain complexes:
\begin{equation}
\begin{tikzcd}
	{C^\bullet} & {C^\bullet[1]} \\
	{C^\bullet} & {C^\bullet[1]}
	\arrow["{T}", from=1-1, to=2-1]
	\arrow["{-T[1]}", from=1-2, to=2-2]
	\arrow["{g}", from=1-1, to=1-2]
	\arrow["{g}", from=2-1, to=2-2]
\end{tikzcd}
\end{equation}
\end{lemma}
\begin{proof}
The chain-map assertion of $g:C^\bullet\to C^\bullet[1]$ in even degrees reads $A'g=gA$ and in odd degrees reads $Ag=gA'$; both identities follow from Lemma \ref{lem:auxiliary_identities}(e). The diagram commutes because $(-T)g=gT$ by Lemma \ref{lem:auxiliary_identities}(e). The chain map $g$ is a chain isomorphism because $g^2:C^\bullet\to C^\bullet[2]=C^\bullet$ is the identity map, where the equal sign reflects the $2$-periodicity of the chain complex $C^\bullet$.
\end{proof}

Lemma \ref{lem:auxiliary_identities}(b)(c)(d) implies a property of the chain map $T$ that holds universally:

\begin{lemma}\label{lem:exact_on_cohomology}
For any $R$-module $M$, the sequence
\begin{equation}
\dots \map[T^i_M]H^i(C_M^\bullet)\map[T^i_M]H^i(C_M^\bullet)\map[T^i_M] \dots
\end{equation}
is exact for $i=0,1$.
\end{lemma}
\begin{proof}
Lemma \ref{lem:auxiliary_identities}(b) implies that $T^2:C^\bullet\to C^\bullet$ is chain homotopic to zero, with the degree-$(-1)$ homotopy $h^i:C^i\to C^{i-1}$ given by $H$ in even degrees and $H'$ in odd degrees; thus $T^2=d^{i-1}h^i+h^{i+1}d^i$. Therefore, the induced map on the cohomology $(T^i_M)^2: H^i(C_M^\bullet)\to H^i(C_M^\bullet)$ is zero.

To show the exactness, it suffices to work with $i=0$ thanks to Lemma \ref{lem:iso_to_shift} tensored with $M$. It remains to prove $\ker T^0_M \subeq \im T^0_M$, where we recall that $T^0_M: \ker A_M/\im A'_M\to \ker A_M/\im A'_M$. Let $\bbar u\in \ker(T^0_M)\subeq \ker A_M/\im A'_M$. Choose a representative $u\in\ker A_M$. In the matrix calculations below, all matrices act on copies of $M$, and we omit the subscript $M$. Then $\bbar u\in \ker(T^0_M)$ implies $Tu=A'v$ for some $v\in M^2$. Note that
\begin{equation}
\begin{bmatrix}
A & O_{2\times 2}
\end{bmatrix}
\begin{bmatrix}
u\\
v
\end{bmatrix}
=
\begin{bmatrix}
T & -A'
\end{bmatrix}
\begin{bmatrix}
u\\
v
\end{bmatrix}=0 \in M^2.
\end{equation}

Consider
\begin{equation}
s=K\begin{bmatrix}
u\\
v
\end{bmatrix}.
\end{equation}

Then by Lemma \ref{lem:auxiliary_identities}(c),
\begin{equation}
As=AK\begin{bmatrix}
u\\
v
\end{bmatrix}=H \begin{bmatrix}
A & O_{2\times 2}
\end{bmatrix}
\begin{bmatrix}
u\\
v
\end{bmatrix}-T
\begin{bmatrix}
T & -A'
\end{bmatrix}
\begin{bmatrix}
u\\
v
\end{bmatrix}=0\in M^2.
\end{equation}

By Lemma \ref{lem:auxiliary_identities}(d),
\begin{equation}
\begin{aligned}
u-Ts&=\parens*{\begin{bmatrix}
 I_{2\times 2} & O_{2\times 2}
\end{bmatrix}
-TK} \begin{bmatrix}
u\\
v
\end{bmatrix}\\
&=
H\begin{bmatrix}
T & -A'
\end{bmatrix}\begin{bmatrix}
u\\
v
\end{bmatrix}
+A'\begin{bmatrix}
O_{2\times 2} & H
\end{bmatrix}\begin{bmatrix}
u\\
v
\end{bmatrix}\\
&=
A'\begin{bmatrix}
O_{2\times 2} & H
\end{bmatrix}\begin{bmatrix}
u\\
v
\end{bmatrix}\in \im(A'_M).
\end{aligned}
\end{equation}

Hence $s\in \ker(A_M)$ and its cohomology class $\bbar s$ satisfies $T^0_M(\bbar s)=\bbar u$, so $\bbar u\in \im(T^0_M)$.
\end{proof}

In particular, for any $R$-module $M$, there is a distinguished subspace $\bbar W_M\subeq H^0(C_M^\bullet)$ with $\bbar W_M=\ker(T^0_M)=\im(T^0_M)$. It is instructive to think about the filtrations associated to the nilpotent operator $T_M^0$ on $H^0(C^\bullet_M)$:
\begin{equation}\label{eq:nilponent_filtration}
0\subeq \bbar W_M \subeq H^0(C_M^\bullet), \qquad H^0(C_M^\bullet)\map[T^0_M] \bbar W_M \map[T^0_M] 0.
\end{equation}

When $M$ is finite-dimensional over $k$, the rank-nullity theorem implies
\begin{equation}\label{eq:middle_dimensionality}
\dim_k H^0(C_M^\bullet) = 2\dim_k \bbar W_M.
\end{equation}

Taking the preimage of \eqref{eq:nilponent_filtration} with respect to the quotient map $\ker A_M \onto \ker A_M/\im A'_M$, we define
\begin{equation}
W^0(M):=\im A'_M \subeq W^1(M):=\ker A_M\cap T_M^{-1}(\im A'_M) \subeq W^2(M):=\ker A_M.
\end{equation}

Then \eqref{eq:middle_dimensionality} implies $\dim_k W^1(M) = (\dim_k W^0(M)+\dim_k W^2(M))/2$. This middle-dimensionality is the only reason why we need Lemma \ref{lem:exact_on_cohomology}.

We summarize the relations among relevant dimensions associated to $C_M^\bullet$ and $T$ for a finite-dimensional $R$-module $M$. Recall $a(M):=\dim_k \ker A_M$ and $b(M):=\dim_k \im A_M$, and we have $a(M)+b(M)=2\dim_k M$ by the rank-nullity theorem.

\begin{lemma}\label{lem:numerical}
For any $R$-module $M$ with $\dim_k M<\infty$, we have
\begin{enumerate}
\item $\dim_k \ker A'_M = \dim_k \ker A_M = \dim_k W^2(M) = a(M)$;
\item $\dim_k \im A'_M = \dim_k \im A_M = \dim_k W^0(M) = b(M)$.
\item $\dim_k H^0(C_M^\bullet) = \dim_k H^1(C_M^\bullet) = a(M) - b(M)$.
\item $\dim_k W^1(M) = \dfrac{a(M)+b(M)}{2}$.
\end{enumerate}
\end{lemma}
\begin{proof}
Identities (a)(b)(c) are immediate from tensoring Lemma \ref{lem:iso_to_shift} with $M$. The identity (d) is restating the middle-dimensionality observation above. Note that $\dim_k H^0(C_M^\bullet)$ is \emph{a priori} $\dim_k \ker A_M - \dim_k \im A'_M$, so it is necessary to observe (b) to express $\dim_k H^0(C_M^\bullet)$ in terms of the rank and the nullity of $A_M$.
\end{proof}

\subsection{A stratification on $V_d$}
In light of Lemma \ref{lem:numerical}, we have a well-defined stratification on $V_d$ by
\begin{equation}\label{eq:stratification_staircase}
V_d=\bigsqcup_{\substack{a+b=2d\\a\geq b\geq 0}} V\sbp{a,b}, \quad V\sbp{a,b}:=\set{M_d\in V_d: a(M_d)=a, b(M_d)=b}.
\end{equation}

Recall that a point $(X_{d+1},Y_{d+1})$ of $V_{d+1}$ is determined by a point $(X_d,Y_d)\in V_d$ and a vector $u:=\begin{bmatrix}
z\\
w
\end{bmatrix}\in \ker A_{M_d}$, according to the rule \eqref{eq:induction_rule}. We then determine the stratum that $(X_{d+1},Y_{d+1})$ belongs to in terms of $M_d$ and $u$.

\begin{lemma}\label{lem:staircase_strata_induction}
In the notation above, suppose $M_d\in V\sbp{a,b}$. Then
\begin{equation}
(X_{d+1},Y_{d+1}) \in \begin{cases}
V\sbp{a+2,b}, & \text{if }u\in W^0(M_d);\\
V\sbp{a+1,b+1}, & \text{if }u\in W^1(M_d)\setminus W^0(M_d);\\
V\sbp{a,b+2}, & \text{if }u\in W^2(M_d)\setminus W^1(M_d).
\end{cases}
\end{equation}
\end{lemma}

Noting that the dimensions of $W^0(M_d),W^1(M_d), W^2(M_d)$ are $b, (a+b)/2, a$ respectively, Lemma \ref{lem:staircase_strata_induction} immediately implies Theorem~\ref{cor:staircase_formula}.

\begin{proof}[Proof of Lemma \ref{lem:staircase_strata_induction}]
It suffices to determine $b(M_{d+1})$, and moreover, we shall use the definition $b(M_{d+1})=\dim_k \im A'_{M_{d+1}}$. The matrix for $A'_{M_{d+1}}$ on $M_{d+1}\oplus M_{d+1}=M_d\oplus k \oplus M_d \oplus k$ is
\begin{equation}
A'_{M_{d+1}} = \begin{bmatrix}
X_{d+1} & Y_{d+1}^2 \\
Y_{d+1} & X_{d+1}
\end{bmatrix}
=
\begin{bmatrix}
X_d & z & Y_d^2 & Y_d w\\
0 & 0 & 0 & 0\\
Y_d & w & X_d & z\\
0 & 0 & 0 & 0
\end{bmatrix}.
\end{equation}

Therefore, the image of $A'_{M_{d+1}}$ actually lies in $M_d\oplus M_d$, and is equal to
\begin{equation}
\begin{aligned}
\im A'_{M_{d+1}} &= \im \begin{bmatrix}
X_d & Y_d^2 \\
Y_d & X_d
\end{bmatrix}
+
\Span_k \set*{\begin{bmatrix}
z\\w
\end{bmatrix},\begin{bmatrix}
Y_d w\\z
\end{bmatrix}}
\\
&= \im A'_{M_d} + \Span_k \set{u, T_{M_d} u}\in M_d\oplus M_d.
\end{aligned}
\end{equation}

We now work mod $\im A'_{M_d}$. Since $u\in \ker A_{M_d}$, we are working in $\ker A_{M_d}/\im A'_{M_d} = H^0(C_{M_d}^\bullet)$. Hence, the increase in the $b$-parameter is given by
\begin{equation}
b(M_{d+1})-b(M_d)=\dim_k \Span_k \set{\bbar u, T^0_{M_d} \bbar u}\subeq H^0(C_{M_d}^\bullet).
\end{equation}

We recall that $T^0_{M_d}$ is a square-zero linear map on $H^0(C_{M_d}^\bullet)$. By a standard fact about nilpotent linear maps (Lemma \ref{lem:cyclic_subspace}), $\dim_k \Span_k \set{\bbar u, T^0_{M_d} \bbar u}$ is the number of nonzero elements in $\set{\bbar u, T^0_{M_d} \bbar u}$. This completes the proof.
\end{proof}

\begin{lemma}\label{lem:cyclic_subspace}
Let $V$ be a finite-dimensional vector space over a field $k$, and $T:V\to V$ be a nilpotent linear map. For any vector $v\in V$, consider the sequence $v,Tv,T^2 v,\dots$, which is eventually zero. Then the nonzero elements in the sequence are linearly independent.
\end{lemma}
\begin{proof}
Assume $T^{n-1}v\neq 0$ and $T^nv=0$. Let $W=\Span_k \set{v,Tv,\dots, T^{n-1}v}$, so $T$ is a nilpotent endomorphism of $W$. We note that $T^{n-1}$ is not the zero map on $W$ because $T^{n-1}v\neq 0$. Hence $\dim W > n-1$, so that $v,Tv,\dots, T^{n-1}v$ must be linearly independent.
\end{proof}

\begin{example}
Using Theorem \ref{cor:staircase_formula}, we compute the motive of $V_d$ for $d\leq 8$ in Table \ref{table:staircase}. The cases with $d\leq 3$ can be directly verified from the definition \eqref{eq:upper_triangular}.
\begin{center}
\begin{table}[h]
\begin{tabular}{|r|l|}
\hline
$d$ & $[V_d]$ \\
\hline
0 & 1\\
1 & 1\\
2 & $\L^2$\\
3 & $3\L^4-2\L^3$\\
4 & $2 \L^8+3 \L^7-5 \L^6+\L^5$\\
5 & $10 \L^{12}-5 \L^{11}-9 \L^{10}+5 \L^9$\\
6 & $5 \L^{18}+21 \L^{17}-30 \L^{16}-9 \L^{15}+15 \L^{14}-\L^{12} $\\
7 & $35 \L^{24}+7 \L^{23}-84 \L^{22}+15 \L^{21}+35 \L^{20}-7 \L^{18}$\\
8 & $14 \L^{32}+112 \L^{31}-112 \L^{30}-162 \L^{29}+113 \L^{28}+70 \L^{27}-7 \L^{26}-28 \L^{25}+\L^{22}$\\
\hline
\end{tabular}

\caption{The motive of $V_d$ for $d\leq 8$}
\label{table:staircase}
\end{table}
\end{center}

We observe from the table that as a polynomial in $\L$, the total coefficient of $[V_d]$ is $1$. Indeed, if we substitute $\L\mapsto 1$ in Theorem \ref{cor:staircase_formula}, we have $\ev{[V_{2d,0}]}{\L\mapsto 1}=1$ and $\ev{[V_{a,b}]}{\L\mapsto 1}=0$ if $(a,b)\neq (2d,0)$, so $\ev{[V_d]}{\L\mapsto 1}=1$. Geometrically speaking, this means that when $k=\C$, the Euler characteristic of $V_d$ in the analytic topology is $1$, which is expected because the obvious $\C^\times$-action $t\cdot (X,Y)=(t^3 X,t^2 Y)$ on $V_d$ has only one fixed point $(X,Y)=(0,0)$.
\end{example}

\begin{remark}[An exponential generating function]
For a quasisymmetric function $Z(t_1,t_2,\dots)$, define
\[ (\mathcal{B}Z)(t)\coloneqq \sum_{d\geq 0} [t_1\cdots t_d]Z(t_1,\dots) \frac{t^d}{d!}.\]
Then the operator $\mathcal{B}$ is multiplicative. A natural generating function that packages $[V_d]$ is thus
\[
(\mathcal{B}\Zhat_R)(t)
:=
\sum_{d\geq 0}
[\Coh_{(1,\dots,d)}(R)]\frac{t^d}{d!}
=
\sum_{d\geq 0}
\frac{[V_d]}{(\L-1)^d\L^{\binom d2}}\frac{t^d}{d!},
\]
where $\Zhat_R(t_1,\dots)$ is viewed as a quasi-symmetric function defined in \S\ref{subsubsec:limit}.

By Theorem~\ref{eq:nest cusp}(d), taking the coefficient of
$t_1\cdots t_d$ gives
\[
 \frac{[V_d]}
 {(\mathbb L-1)^d\mathbb L^{\binom d2}}
 =
 e_0^T\widehat H_{R,1}^{d}\mathbf 1,
\]
the sum of entries in the $0$-th row of $H_{R,1}^d$.
Consequently,
\[
 (\mathcal{B}\Zhat_R)(t)=\sum_{d\geq 0}
 \frac{[V_d]}
 {(\mathbb L-1)^d\mathbb L^{\binom d2}}
 \frac{t^d}{d!}
 =
 e_0^{T}
 \exp\!\left(t\widehat H_{R,1}\right)\mathbf 1,
\]
getting a generating function of $[V_d]$ that can be checked against the recursions in Theorem~\ref{cor:staircase_formula}.
\label{rmk:borel}
\end{remark}
\bibliographystyle{alpha}
\bibliography{bib}

\end{document}